\date{}
\documentclass [11pt]{article}
\usepackage{amsfonts,amsmath,amssymb,graphicx}

\title{The early evolution of the $H$-free process}
\author{
Tom Bohman \thanks{Department of Mathematical Sciences,
Carnegie Mellon University, Pittsburgh, USA. Supported in part by
NSF grant DMS-0701183. E-mail: \texttt{tbohman@math.cmu.edu}}
\and
Peter Keevash \thanks{School of Mathematical Sciences,
Queen Mary, University of London, Mile End Road, London E1 4NS, UK.
Email: p.keevash@qmul.ac.uk.
Research supported in part by NSF grant DMS-0555755.}
}

\newtheorem{theo}{Theorem}[section] 

\newtheorem{lemma}[theo]{Lemma}
\newtheorem{coro}[theo]{Corollary}
\newtheorem{conj}[theo]{Conjecture}

\newcommand{\mc}[1]{\mathcal{#1}}
\newcommand{\mb}[1]{\mathbb{#1}}
\newcommand{\nib}[1]{\noindent {\bf #1}}

\newcommand{\sm}{\setminus}
\newcommand{\ov}{\overline}

\newcommand{\eps}{\epsilon}

\newcommand{\sub}{\subseteq}
\newcommand{\subn}{\subsetneq}

\newcommand{\GG}{\Gamma}
\newcommand{\TT}{{\mathcal T}}

\newcommand{\aA}{\alpha}
\newcommand{\bB}{\beta}
\newcommand{\gG}{\gamma}

\newcommand{\tT}{\theta}

\def\COMMENT#1{}

\def\qed{\hfill $\Box$}

\begin{document}
\maketitle

\begin{abstract}
The $H$-free process, for some fixed graph $H$, is the random graph process
defined by starting with an empty graph on $n$ vertices 
and then adding edges one at a time,
chosen uniformly at random subject to the constraint that no 
$H$ subgraph is formed.
Let $G$ be the random maximal $H$-free graph
obtained at the end of the process. When $H$ is strictly $2$-balanced,
we show that for some $c>0$, with high probability as $n \to \infty$,
the minimum degree in $G$ is at least $cn^{1-(v_H-2)/(e_H-1)}(\log n)^{1/(e_H-1)}$.
This gives new lower bounds for the Tur\'an numbers of certain bipartite graphs,
such as the complete bipartite graphs $K_{r,r}$ with $r \ge 5$.
When $H$ is a complete graph $K_s$ with $s \ge 5$  we show that for some $C>0$,
with high probability the independence number of $G$ is at most
$Cn^{2/(s+1)}(\log n)^{1-1/(e_H-1)}$. This gives new lower bounds for Ramsey
numbers $ R(s,t)$  for fixed $s \ge 5$ and $t$ large.
We also obtain new bounds for the independence number of $G$ for other graphs $H$,
including the case when $H$ is a cycle. Our proofs use the differential equations
method for random graph processes to analyse the evolution of the process, and 
give further
information about the structure of the graphs obtained, including asymptotic
formulae for a broad class of subgraph extension variables.
\end{abstract}

\section{Introduction}

Random graph processes provide a natural context for modeling a complex network that
evolves over time.  While there has been considerable recent interest in
using such processes to model networks that arise in applications
(see \cite{D} and the references therein), random graphs have long been an important
component in the construction of sophisticated 
combinatorial objects (see \cite{AS}). 
In the classical Erd\H{o}s-R\'enyi random graph model $G(n,p)$ each
pair of vertices appears as an edge with probability $p=p(n)$ and these choices
are mutually independent.  The closely related random graph $ G(n,i) $ is chosen
uniformly at random from the collection of all 
graphs with \(n\) vertices and \(i\) edges.
These models are well understood, but
distributions on graphs given by random processes in which there is significant dependence 
among the choices made in different rounds are typically much more difficult to analyse.
For many such processes even the most basic quantities, such as the number
of edges in the final graph, are not known (see \cite{G}, for example).

In this paper we analyse a significant portion of the
initial evolution of the $H$-free process, for some fixed graph $H$,
defined by starting with an empty graph on $n$ vertices and then adding edges one at a time,
chosen uniformly at random subject to the constraint that no $H$ subgraph is formed.
More formally, we begin with the graph on $n$ vertices with no edges,
which we denote $G(0)$. Now suppose $i>0$ and we have some graph $G(i-1)$. We say that
a pair $uv$ of vertices is {\em open} in $G(i-1)$ if $uv$ is not an edge of $G(i-1)$
and $G(i-1) \cup \{uv\}$ does not contain $H$ as a subgraph. We choose $uv$ uniformly
at random among all open pairs in $G(i-1)$ and then $G(i)$ is obtained from $G(i-1)$
by adding the edge $e_i = uv$. The process terminates when there are no open pairs,
with some graph $G$ on $n$ vertices that is a maximal $H$-free graph.
Beside being of interest in its own right, our analysis of this process produces
new results in Ramsey theory and the theory of Tur\'an problems.

Erd\H{o}s, Suen and Winkler \cite{ESW} suggested this process as a means
to generate an interesting probability distribution on the collection of 
maximal $H$-free graphs, or more generally maximal graphs with any fixed graph property.%
\footnote{Bollob\'as (personal communication) informs us that such processes 
were considered earlier, if not in print.}
They obtained results on the triangle-free process and the bipartite process, 
using a differential equations method that had been previously applied by 
Ruci\'nski and Wormald \cite{RW} to analyse the `maximum degree $d$' process.
Another motivation for their work was that
their analysis of the triangle-free process led to the best lower bound on the
Ramsey number $R(3,t)$ known at that time.

Ramsey theory encompasses a variety of results expressing the informal principle
that all large systems have some structure. It is a source of many challenging
unsolved combinatorial problems and has applications throughout mathematics.
We refer the reader to \cite{GRS} for an introduction to the subject.
The Ramsey number $R(s,t)$ is the least number $n$ such that any graph on $n$ vertices
contains a complete graph with $s$ vertices or an independent set with $t$ vertices.
In general, very little is known about these numbers, even approximately.
The upper bound $R(3,t) = O(t^2/\log t)$ was obtained by Ajtai, Koml\'os and
Szemer\'edi \cite{AKS}, but for many years the best known lower bound, due to Erd\H{o}s \cite{E},
was $\Omega(t^2/\log^2 t)$. Spencer conjectured that the triangle-free
process is likely to produce a graph that establishes a good lower bound on \( R(3,t) \) for
$t$ large; the idea being that the triangle-free process admits enough random edges to bring the
independence number close to the smallest possible for a triangle-free graph.
Finally, Kim \cite{K} determined the order of magnitude, showing
that $R(3,t) = \Theta(t^2/\log t)$. His proof made use of a semi-random construction that
is motivated (even guided) by the triangle-free process, but the
question remained open as to whether the triangle-free process itself gives such
a good construction. This was answered by Bohman \cite{B}, who showed that
with high probability, the graph produced by the triangle-free process 
has independence number bounded above by $O(n^{1/2} \log^{1/2} n)$ 
and minimum degree bounded below by $\Omega(n^{1/2}\log^{1/2} n)$.
He went on to analyse the $K_4$-free process, improving the best known lower bound
on $R(4,t)$ to $R(4,t) > \Omega(t^{5/2}/\log^2 t)$.

The general $H$-free process was independently studied by Osthus and Taraz \cite{OT}
and by Bollob\'as and Riordan \cite{BR}. Say that a graph $H$ is
{\em strictly $2$-balanced} if the number of vertices $v_H$ and edges $e_H$ in $H$
are both at least $3$ and
$$\frac{e_H-1}{v_H-2} > \frac{e_K-1}{v_K-2}$$
for all proper subgraphs $K$ of $H$ with $v_K \ge 3$.
Osthus and Taraz showed that if $H$ is strictly $2$-balanced then for some $c,C>0$
with high probability, for the $H$-free process $G$ has
average degree at least $cn^{1-(v_H-2)/(e_H-1)}$
and maximum degree at most $Cn^{1-(v_H-2)/(e_H-1)}(\log n)^{1/(\Delta(H)-1)}$.
(In fact they proved the average degree bound under a similar but weaker condition on $H$.)
Wolfovitz \cite{Wz} showed that if $H$ is strictly $2$-balanced and regular then
the expected number of edges in $G$ is at least $cn^{2-(v_H-2)/(e_H-1)}(\log \log n)^{1/(e_H-1)}$.
An immediate consequence is an improved lower bound for Tur\'an numbers,
which leads us to another motivation for studying the $H$-free process.

The Tur\'an number $\mbox{ex}(n,H)$ is the maximum possible number of edges in a graph
on $n$ vertices that does not contain an $H$ subgraph. More generally, the theory
of Tur\'an problems concerns the study of combinatorial structures that have
maximum size subject to not containing some fixed structure. We refer the
reader to \cite{F} for a survey of this subject. Tur\'an \cite{T} determined the value
of $\mbox{ex}(n,H)$ when $H=K_r$ is complete: the unique largest graph on $n$ vertices with
no $K_r$ subgraph is complete $(r-1)$-partite with part sizes as equal as possible.
For general $H$, the Erd\H{o}s-Stone-Simonovits theorem \cite{ESt, ESi} gives
the estimate $\mbox{ex}(n,H)=\mbox{ex}(n,K_r)+o(n^2)$, where $r=\chi(H)$ is the
chromatic number of $H$. This gives an asymptotic formula for the Tur\'an number
when $H$ is not bipartite. However, when $H$ is bipartite it is an open problem in
general to determine even the order of magnitude of $\mbox{ex}(n,H)$. For example,
when $H=K_{r,r}$ is complete bipartite with $r \ge 5$, for many years the best known
lower bound was $\mbox{ex}(n,K_{r,r}) = \Omega(n^{2-2/(r+1)})$, a result of
Erd\H{o}s and Spencer \cite{ESp} proved via a simple application of the probabilistic method.
Wolfovitz's analysis of the $H$-free process improved this to
$\mbox{ex}(n,K_{r,r}) = \Omega(n^{2-2/(r+1)}(\log \log n)^{1/(r^2-1)})$.

\subsection{Results I: Ramsey and Tur\'an bounds}

In this paper we extend the methods from \cite{B} to an analysis of
the $H$-free process when $H$ is strictly $2$-balanced, leading to new lower
bounds for Ramsey and Tur\'an numbers. We also investigate other properties of
the process, viewing it as a model of interest in its own right, and give certain
extension counting formulae that address a question of Spencer.  In particular, we 
show that the graph produced by the \(H\)-free process is very similar to the 
corresponding random graph \( G(n,i) \) with respect to
small subgraph counts, with the exception 
that the \(H\)-free process produces no copies of graphs containing \(H\).
We begin with the Tur\'an and Ramsey results.

Our first theorem gives a new lower bound for the number of edges in $G$.
In fact we have a new lower bound for the minimum
degree, and it holds with high probability, not just in expectation.
An immediate consequence is a lower bound for the Tur\'an number $\mbox{ex}(n,H)$.

\begin{theo} \label{turan}
Suppose that $H$ is a strictly $2$-balanced graph with $v_H$ vertices and $e_H$ edges.
Then for some $c>0$ with high probability the minimum degree in the final graph of the
$H$-free process is at least $cn^{1-(v_H-2)/(e_H-1)}(\log n)^{1/(e_H-1)}$. In particular,
the Tur\'an number satisfies
$$\mbox{ex}(n,H) = \Omega \left(n^{2-(v_H-2)/(e_H-1)}(\log n)^{1/(e_H-1)}\right).$$
\end{theo}
\noindent
Note that it follows immediately from Theorem~\ref{turan} that we have
$$ \mbox{ex}\left(n, K_{r,r} \right) = \Omega \left(n^{2-2/(r+1)}(\log n)^{1/(r^2-1)} \right).$$
For general complete bipartite graphs $K_{r,s}$ with $r \le s$,
the `Zarankiewicz problem' of estimating $\mbox{ex}\left(n,K_{r,s}\right)$ 
is a subject of special interest in extremal graph theory. A general upper bound
of order $n^{2-1/r}$ was given by K\"ov\'ari, S\'os and Tur\'an \cite{KST}.
The only known asymptotic results are 
$\mbox{ex}\left(n,K_{2,r}\right) \sim \frac{1}{2}(r-1)^{1/2}n^{3/2}$
(see \cite{F2}) and
$\mbox{ex}\left(n,K_{3,3}\right) \sim \frac{1}{2}n^{5/3}$
(see \cite{Br} and \cite{F3}).
Note that the lower bound construction for $K_{3,3}$ also gives the 
best known lower bound for $K_{4,4}$.
The only other case when the upper bound is known to be of the
correct order of magnitude is when $s>(r-1)!$ (see \cite{ARS}).
The known constructions are based on algebraic and geometric structures that may
not exist for other values of the parameters $r$ and $s$.
However, it is widely believed that \( \mbox{ex}\left(n, K_{r,s}\right) \) 
for general $r \le s$ is on the order of \( n^{2-1/r} \).

For Ramsey numbers, we obtain the following new lower bounds.
\begin{theo} \label{ramsey}
For fixed $s \ge 5$ and $t \to \infty$, the Ramsey number satisfies
$$R(s,t) = \Omega \left(t^{\frac{s+1}{2}}(\log t)^{\frac{1}{s-2}-\frac{s+1}{2}} \right).$$
\end{theo}
\noindent
The previously best known lower bound on \( R(s,t) \) when \(s\) is fixed and \(t\) is large
was $R(s,t) = \Omega \left( (t/\log t)^{\frac{s+1}{2}} \right)$,
established by Spencer \cite{Sp1} using the Lov\'asz Local Lemma.
Theorem \ref{ramsey} improves this by a multiplicative factor of \( (\log t)^{1/(s-2)} \).
There is no particular reason to believe that our lower bound is anywhere near optimal,
since the best known general upper bound is essentially $t^{s-1}$
(up to a polylogarithmic factor in $t$).  On the other hand, as Theorem~\ref{ramsey} 
can be viewed as the natural generalisation of the construction that gives the correct order 
of magnitude for \( R(3,t) \), it would be interesting to see a significant improvement 
on the bound in Theorem~\ref{ramsey} for \( s \ge 4 \).

We also obtain new lower bounds for cycle-complete Ramsey numbers. Given graphs $H_1$, $H_2$,
the graph Ramsey number $R(H_1,H_2)$ is the least number $n$ such that for any $2$-colouring of the
edges of $K_n$ there is a monochromatic copy of $H_1$ or $H_2$. Note that $R(C_\ell,K_t) \ge n$
if and only if there is a $C_\ell$-free graph on $n$ vertices with no independent set of size $t$.
We prove the following bound.

\begin{theo} \label{cycle-ramsey}
For fixed $\ell \ge 4$ and $t \to \infty$ the cycle-complete Ramsey number satisfies
$$R(C_\ell,K_t) = \Omega \left( (t/\log t)^{\frac{\ell-1}{\ell-2}} \right).$$
\end{theo}
\noindent
Again this is quite far from the best known upper bounds (see \cite{CLRZ,LZ,Su}). 
For example, Erd\H{o}s \cite{E2} conjectured that $R(C_4,K_t) = O(t^{2-\eps})$ 
for some absolute constant $\eps>0$, but this is still open.

In fact, we establish more general properties of the $H$-free process from
which these theorems follow. In order to show that the process continues to
run for a certain number of steps, we will establish asymptotic formulae
for various graph parameters at any given time in the process, including
the degree of any vertex, but also more general extension parameters.
To state these formulae we need some terminology and notation.

\subsection{Terminology and notation I} \label{notate1}

We write $[n]=\{1,\cdots, n\}$ for the vertex set of the process.
At step $i$ of the process let $E(i)$ be the edges of the graph $G(i)$,
let $O(i)$ be the pairs of vertices that are open (as defined above), and
let $C(i)$ be the pairs of vertices that are neither edges nor open,
which we refer to as {\em closed}.

We fix some strictly $2$-balanced graph $H$ throughout the paper and
write $$p = n^{-\frac{v_H-2}{e_H-1}}.$$
For any graph $\GG$ we write $V_\GG$ for the vertex set of $\GG$, $E_\GG$ for the
edge set of $\GG$, $v_\GG = |V_\GG|$ and $e_\GG=|E_\GG|$.
For $A \sub V_\GG$ we write
$$S_\GG = p^{e_\GG} n^{v_\GG}  \ \ \ \ \  \text{ and }  \ \ \ \ \  S_{A,\GG} = p^{e_\GG-e_{\GG[A]}} n^{v_\GG-|A|}.$$
We say that such a pair $(A,\GG)$ is 
{\em strictly balanced} if $S_{A,\GG[B]}>S_{A,\GG}$ for every $A \subn B \subn V_\GG$
and {\em strictly dense} if $S_{A,\GG[B]}>1$ for every $A \subn B \sub V_\GG$.  

A key element of our analysis of the \(H\)-free process is closely tracking 
the number of extensions from  fixed sets of vertices to fixed subgraphs of \( G(i) \).  
Intuitively, the graph \( G(i) \) produced by the
\( H \)-free process should be roughly equal to the random graph \( G(n,i) \), the graph chosen 
uniformly at random from the collection of graphs with \(n\) vertices and \(i\) edges,
up until the number of copies of \(H\) in \( G(n,i) \) is roughly equal to the number of edges.
This occurs when \( i \) is roughly \( pn^2 \), with \(p\) as defined above. We expect 
the more interesting part of the evolution of the \(H\)-free process to be 
at and beyond this range of \(i\). Considering $G(n,p)$, which is very 
similar to $G(n,i)$ here, we note that
\( S_\GG \) is roughly the expected number of labeled copies of $\GG$, 
and \( S_{A,\GG} \) is roughly the expected number of labeled extensions to $\GG$
from a fixed set of vertices playing the role of $A$.
Thus we can think of these quantities as anticipated {\em scalings} by which
we should measure the same parameters in the $H$-free process.

In order to track extensions, we track all `open routes' to such extensions.
Suppose $\GG$ is a graph and $J$ is a spanning subgraph of $\GG$.
Suppose also that $A \sub V_\GG$ is an independent set
in $\GG$ and $\phi:A \to [n]$ is an injective mapping.
We define the {\em extension variables} $X_{\phi,J,\GG}(i)$ to be
the number of injective maps $f:V_\GG \to [n]$ such that
\begin{itemize}
\item[(i)] $f(e) \in O(i)$ for every $e \in E_\GG \setminus E_J$,
\item[(ii)] $f(e) \in E(i)$ for every $e \in E_J$, and
\item[(iii)] $f$ restricts to $\phi$ on $A$.
\end{itemize}
We say that the random variable $X_{\phi,J,\GG}(i)$ is {\em trackable}
if one of the following two conditions holds:
\begin{itemize}
\item[(a)] $ (A,\GG) $ is strictly dense and \( \GG \) does not contain \(H\) as a subgraph, or
\item[(b)] $ S_{A, \GG} =1 $, $(A,\GG)$ is strictly balanced, \( E_J \subsetneq E_\GG \),
and $H$ is not a subgraph of the graph $\GG'$
obtained from $\GG$ by adding the edges $ab$ for all $a,b \in A$ with $\phi(a)\phi(b) \in E(i)$.
\end{itemize}
It follows easily from the definitions 
that for any trackable extension variable \( X_{ \phi, J, \GG}(i) \) the
pair \( (A,J) \) is strictly dense.
Note further that condition (b) includes the case where
$\GG = H \sm ab$ for some $ab \in E_H$, $e_J \le e_H-2$,
$A = \{a,b\}$ and $\phi(ab) \notin E(i)$.  These extensions comprise the set of {\em open routes}
to a copy of \(H\) less an edge, where \(\phi(ab)\) plays the role 
of the missing edge.  As the appearance of such an extension is the mechanism whereby 
the pair \(\phi(ab)\) becomes closed, these particular extension variables play a
central role in our analysis of the \(H\)-free process.

We fix constants $V, W, \eps, \mu$ throughout the paper which satisfy
$0 < \mu \ll \eps \ll 1/W \ll 1/V \ll 1/e_H$. 
(The notation $0 < \aA \ll \bB$ means that there is an increasing function
$f(x)$ so that the following argument is valid for $0 < \aA < f(\bB)$.)
We introduce a continuous time
variable $t$, using the scaling $t = t(i) = i/s$ with $s=pn^2$,
and analyse the process up to time $t_{\max} = \mu (\log n)^{1/(e_H-1)}$,
which corresponds to
$$m = \mu (\log n)^{1/(e_H-1)} pn^2$$
edges. Let $\TT$ be the set of all triples $(A,J,\GG)$ where
$J$ is a spanning subgraph of a graph $\GG$ with $v_\GG, e_\GG < V$, 
$A$ is an independent set in $\GG$, and the variables $X_{\phi,J,\GG}(0)$
are trackable. Write $aut(H)$ for the number of automorphisms of $H$ and define
$$q(t) = e^{-2e_H aut(H)^{-1} (2t)^{e_H-1}},
\quad P(t) = W(t^{e_H-1} + t), \quad e(t) = e^{P(t)}-1
\quad \mbox{ and } \quad s_e = n^{1/2e_H-\eps}.$$
We also define $\gG(t)$ to be any smooth increasing function such that
$\gG(t) = 40Ve^{40V} t $ for $0 \le t \le 40V/W$,
$ \gG'(t) > 20V $ for $ 40V/W < t \le 1/(50V) $,
and $\gG(t) < 1/2$, $\gG'(t) < W$ for all $t \ge 0$. 
Then we set $\tT(t) = 1/2 + \gG(t)$, 
so that $1/2 \le \tT(t) < 1$ for all $t \ge 0$.

\subsection{Results II: The \(H\)-free process}

Our first main theorem gives asymptotic formulae for trackable extension variables
throughout the process.
\begin{theo} \label{track}
With high probability, for every $i \le m$ and trackable extension variable
$X_{\phi,J,\GG}(i)$ corresponding to a triple in $\TT$, we have
$$X_{\phi,J,\GG}(i) = (1\pm e(t)/s_e)(x_{A,J,\GG}(t) \pm 1/s_e) S_{A,J},$$
where $$x_{A,J,\GG}(t) = (2t)^{e_J}q(t)^{e_\GG-e_J}.$$
\end{theo}
\noindent
(For this theorem to be useful we choose $\eps < \eps(V)$ sufficiently small and
then $\mu < \mu(\eps)$ sufficiently small so that $e(t)$ and $q(t)^{-V}$
are both at most $n^{\eps}$ for $t \le t_{\max}$.)  
Note, for example, that there is a trackable extension variable describing the
number of common neighbours of a set of size $d$ 
whenever $p^dn>1$, so we have the following corollary.
\begin{coro} \label{neighbours}
With high probability, for every $d$ with $p^dn>1$, set $A$ of $d$ vertices
and $i \le m$, the number of common neighbours of $A$ in $G(i)$
is $(1+o(1))(2i/n^2)^d n$.
\end{coro}
\noindent
A remarkable consequence of Theorem~\ref{track} is that the graph \( G(i) \) for 
\( i \le m \) is similar to the uniform random graph \( G(n,i) \) with respect 
to small subgraph counts, with the notable exception that there are no copies of graphs containing
\( H\) in \( G(i)\). The possibility of this intriguing behavior was first suggested 
by Joel Spencer. The following theorem gives the correct asymptotic counts for labelled
copies of a graph $\GG$ in the `subcritical' case (i) and the `supercritical' case (ii).
For the sake of brevity we just establish existence of a copy in the `critical' case (iii),
although our discussion in Section \ref{sec:count} points the way towards better
results in this case.

\begin{theo} \label{count}
Suppose $\GG$ is an $H$-free graph and write $X_\GG(i)$ for the
number of labelled copies of $\GG$ in $G(i)$. Then with high probability
\begin{itemize}
\item[(i)] If there exists $B \sub V_\GG$ with $S_{\GG[B]}<1$ then $X_\GG(m)=0$.
\item[(ii)] If $S_{\GG[B]}>1$ for all non-empty $B \sub V_\GG$ then
$X_\GG(i) \sim (2i/n^2)^{e_\GG}n^{v_\GG}$.
\item[(iii)] If $S_{\GG[B]} \ge 1$ for all $B \sub V_\GG$ then $X_\GG(m)>0$.
\end{itemize}
\end{theo}

While Theorem~\ref{track} alone is enough to establish the Tur\'an bounds stated above, 
our results on the Ramsey numbers require an upper bound on the 
independence number of \( G(m) \). Theorem \ref{ramsey} follows easily
from \ref{k-indep} below. This in turn follows from the following more
general result for $s \ge 6$. (Then we will need to modify the proof slightly
to deal with the case $s=5$.)

\begin{theo} \label{balanced-indep}
Suppose that $H$ is strictly $2$-balanced and that
for any two edges $uv$, $xy$ of $H$ and $\{x,y\} \subn B \subn V_H$
we have $S_{B,H \sm uv} < 1$. Then there is $C>0$ such that with high probability
the final graph of the $H$-free process has independence number at most
$Cn^{(v_H-2)/(e_H-1)}(\log n)^{1-1/(e_H-1)}$.
\end{theo}

\begin{theo} \label{k-indep}
For any $s \ge 5$ there is $C>0$ such that with high probability the final graph of the
$K_s$-free process has independence number at most
$Cn^{\frac{2}{s+1}}(\log n)^{1-\left(\binom{s}{2}-1\right)^{-1}}$.
\end{theo}

Alon, Ben-Shimon and Krivelevich \cite{ABK} recently proposed a construction that takes a
nearly regular \(K_s\)-free graph \(G\) and produces a regular \(K_s\)-free graph with
roughly the same independence number as the original graph.  It follows from
Corollary~\ref{neighbours} that the graph produced after $m$ steps of the \(K_s\)-free process is
a suitable input for this construction. This suggests that the bound on \( R(s,t) \) given 
in Theorem~\ref{ramsey} can be achieved by a regular graph. (A formal proof would need to provide
some details missing from the sketch given in \cite{ABK}.)

We also obtain the following bound when $H$ is a cycle, which implies
Theorem \ref{cycle-ramsey}.

\begin{theo} \label{cycle-indep}
For any $\ell \ge 3$ there is $C>0$ such that with high probability
the final graph of the $C_\ell$-free process has independence number at most
$C(n\log n)^{(\ell-2)/(\ell-1)}$.
\end{theo}

\subsection{Organisation of the paper}

In the next section we give a heuristic explanation for
the differential equations leading to the formulae in Theorem \ref{track}.
In Section 3 we develop some theory of strictly $2$-balanced graphs and balanced extensions. 
Over the following three sections we collect various properties
that hold with high probability on the `good' event at a given time that the process
has followed the trajectory of the differential equations so far.
Section 4 contains various union bound arguments, Section 5 gives upper bounds
on the extension variables and Section 6 provides a means to approximate the
number of pairs that become closed when some particular pair is added as an edge.
In Section 7 we formulate our framework for showing that the process follows
the differential equations, which is based to some extent on that given by
Wormald \cite{W}, but also incorporates martingale estimates from \cite{B}.
Section 8 concerns trackable random variables: we obtain bounds
on the one-step changes of trackable random variables sufficient 
to apply the differential equations method. Then we apply the differential equation
method in Section 9 to prove Theorem \ref{track}, from which 
Theorem \ref{turan} immediately follows. We also apply Theorem \ref{track}
to prove Theorem \ref{count} in Section 10.
Next we turn our attention to the independence number.
In Section 11 we formulate a general property, which we call
`smooth independence', and bound the independence number under
the assumption that $H$ has this property. Then in Section 12 we show
that cycles and complete graphs $K_s$, $s \ge 5$ have smooth independence,
from which Theorems \ref{cycle-indep} and \ref{ramsey} follow.
We also prove Theorem \ref{balanced-indep} in this section.
The final section contains some concluding remarks.

\subsection{Terminology and notation, II}

We write $\mc{G}_j$ for the {\em good} event that 
for every $0 \le i \le j$ and trackable extension variable
$X_{\phi,J,\GG}(i)$ corresponding to a triple in $\TT$, we have
$$X_{\phi,J,\GG}(i) = (1\pm e(t)/s_e)(x_{A,J,\GG}(t) \pm \tT(t)/s_e) S_{A,J}.$$
Note that this implies the formulae in the statement of
Theorem \ref{track}, since $\tT(t) < 1$ for all $t \ge 0$.

When we count extensions it is convenient to work with labeled graphs,
and we will often write $uv$ for the ordered pair $(u,v)$ as well as the edge $\{u,v\}$.
The prime symbol $'$ is occasionally used to denote differentiation with respect
to the time variable $t$: this will be clear from the context.

Statements containing the symbols $\pm$ and/or $\mp$
are shorthand for two separate statements: one with every $\pm$ replaced by $+$ and every $\mp$ by $-$,
the other with $\pm$ replaced by $-$ and $\mp$ by $+$. We also use the notation $a = b \pm c$
to mean $b-c<a<b+c$. Where there is possibility for confusion we label the
symbols as $\pm_1$ and $\pm_2$, e.g. $a^{\pm_1 \pm_2} = b^{\pm_1} \pm c^{\mp_2}$
is shorthand for $4$ separate statements, one of which is $a^{++} = b^+ \pm c^-$.

The parameter $n$ will always be sufficiently large compared to all other
parameters, and we use the phrase `with high probability' to refer to an
event that has probability $1-o_n(1)$, i.e. the probability tends to $1$
as $n$ tends to infinity. In fact we can arrange that our high probability
events fail with probability at most $\exp(-n^\eps)$.

We say that a graph $W$ is a {\em join} of two graphs $W_1$ and $W_2$ if it
has subgraphs $J_1$ isomorphic to $W_1$ and $J_2$ isomorphic to $W_2$
such that $V_W = V_{J_1} \cup V_{J_2}$ and $E_W = E_{J_1} \cup E_{J_2}$.
For convenient notation we use names for vertices in $J_1$ interchangeably
with their corresponding vertices in $W_1$, and similarly for $J_2$ and $W_2$.

If $X$ is a set and $k$ is a non-negative integer then we write $\binom{X}{k}$
for the set of subsets of $X$ of size $k$.

We will not often refer explicitly to the underlying probability space for the $H$-free process,
but we note here the following natural construction. Let $\Omega = \Omega_n$ be the
set of all maximal sequences in $\binom{[n]}{2}$ with distinct entries and the property that
each initial sequence gives an $H$-free graph on vertex set $[n]$.
We stress that our measure is not uniform: it is the measure given by
the uniform random choice at each step. We always work with the natural filtration
$\mc{F}_0 \sub \mc{F}_1 \sub \dots$ given by the process.
Two elements $x,y$ of $\Omega$ are in the same atom (i.e. part of the generating partition)
of $\mc{F}_j$ exactly when the first $j$ entries of $x$ and $y$ agree.

\section{Trajectory equations}

We start by giving a heuristic explanation of the equations describing the
evolution of the $H$-free process. We will then prove the validity of these equations
in subsequent sections.
Recall that $G(i)$ denotes the graph on $[n]$ obtained after
$i$ steps of the $H$-free process: its edge set $E(i)$ contains $i$ edges.
We partition the non-edges $\binom{[n]}{2} \sm E(i)$ into two sets $O(i)$ and $C(i)$,
which we call {\em open} pairs and {\em closed} pairs, respectively.
We say that a pair $uv$ is open if $G(i) \cup uv$ does not contain a copy of $H$,
i.e. $uv$ is a possible choice for the next edge in the process.

\begin{quote}
\nib{Notation.}
We consider the following random variables.
Suppose $\GG$ is a graph and $J$ is a spanning subgraph of $\GG$ (i.e. $V_J=V_\GG$).
Suppose also that $A \sub V_J$ is an independent set (i.e. does not span any edges)
in $\GG$ and $\phi:A \to [n]$ is an injective mapping.
{\em Throughout this paper we assume that $\GG, J, A, \phi$ satisfy these conditions,
even if this is not explicitly stated.}
We define the {\em extension set} $\Xi_{\phi,J,\GG}(i)$ to be
the set of injective maps $f:V_\GG \to [n]$ such that
(i) $f(e) \in O(i)$ for every $e \in E_\GG \setminus E_J$,
(ii) $f(e) \in E(i)$ for every $e \in E_J$, and
(iii) $f$ restricts to $\phi$ on $A$.
Then we define the {\em extension variables} by
$X_{\phi,J,\GG}(i) = |\Xi_{\phi,J,\GG}(i)|$.
In words, we are counting labeled copies (not necessarily induced)
of a graph $J$ in $G(i)$ that extend a particular embedding $\phi:A \to [n]$,
with the extra condition that some extra pairs (i.e. the edges of $\GG \sm J$) are open.
Actually we will be interested in the number of copies up to isomorphism,
but the equations for labeled copies are easier to work with.
\end{quote}

\begin{quote}
\nib{Examples.}
One special case of this definition is the number of labeled copies of a graph $\GG$ in $G(i)$:
this can be written as $X_{\phi_0,\GG,\GG}(i)$,
where we write $\phi_0$ for the unique function $\phi_0: \emptyset \to [n]$.
To count edges and open pairs with this notation we write $e$ and $\ov{e}$
for the two graphs on two vertices, say $\{a,b\}$, with one edge and no edges respectively.
Then $X_{\phi_0,\ov{e},e}(i) = 2|O(i)|$ and $X_{\phi_0,e,e}(i) = 2|E(i)|$.
We can also express the degree $d_{G(i)}(v)$ of a vertex $v$ in $G(i)$
as $X_{\phi_v,e,e}(i)$, where again $e$ is the edge $ab$ and we write $\phi_v$
for the function $\phi: \{a\} \to [n]$ defined by $\phi(a)=v$.
\end{quote}

We write $Q(i)=2|O(i)|$ for the number of ordered pairs that are open.
For an ordered pair $uv \in O(i)$,
write $C_{uv}(i)$ for the set of ordered pairs $xy \in O(i)$ that would become closed,
i.e. belong to $C(i+1)$, if at time $i+1$ the process chooses $uv$ as the edge $e_{i+1}$.
By the definition of $C(i+1)$ this means that adding $uv$ and $xy$ to $G(i)$ would create a copy of $H$.
Another way to say this is that there is a subgraph $J$ obtained by deleting
two edges $ab$ and $cd$ from $H$ and an injective map
$f:V_H \to [n]$ such that $f(a)=u$, $f(b)=v$, $f(c)=x$, $f(d)=y$
and $f(e) \in E(i)$ for every edge of $J$.
We have $f \in \Xi_{\phi_T,J_T,\GG_T}(i)$, where
given such a quadruple $T=(a,b,c,d)$, we write $\GG_T = H \sm ab$, $J_T = H \sm \{ab,cd\}$
and define $\phi_T$ by $\phi_T(a)=u$ and $\phi_T(b)=v$.
In principle there could be many embeddings $f$ giving the same pair $xy$,
but we will show in Lemma \ref{closure} that this is very unlikely:
for most $xy \in C_{uv}(i)$ there will be exactly one such embedding $f$,
up to an automorphism of $H$.
We will see that $C_{uv}(i) \sim aut(H)^{-1} \sum_T X_{\phi_T,J_T,\GG_T}(i)$,
where the sum is over quadruples $T=(a,b,c,d)$ such that
$ab$ and $cd$ are distinct (but not necessarily disjoint) edges of $H$.

To approximate the extension variables we introduce a continuous time variable $t$,
using the scaling $t = t(i) = i/s$ with $s=pn^2$, where we recall that $p = n^{-(v_H-2)/(e_H-1)}$.
We noted above that this is the point at which the number of copies of $H$ in the random
graph $G(n,s)$ is comparable to the number of edges $s$, so it is natural to anticipate
the interesting behaviour to occur at this scale.
We analyse the process up to time $t_{\max} = \mu (\log n)^{1/(e_H-1)}$,
for some small constant $\mu>0$,
which corresponds to $m = \mu (\log n)^{1/(e_H-1)} pn^2$ edges.
For the variable $X_{\phi,J,\GG}(i)$ with $\phi:A \to [n]$
we use the scaling $S_{A,J} = p^{e_J} n^{v_J-|A|}$.
Again, we noted above that the count of these extensions in $G(n,s)$
suggests the use of this scaling.
Our eventual aim is to prove that with high probability, for every $i \le m$
and for every trackable extension variable $X_{\phi,J,\GG}(i)$ corresponding to a triple in $\TT$, 
we have the asymptotic formula
$$X_{\phi,J,\GG}(i) = (1\pm e(t)/s_e)(x_{A,J,\GG}(t) \pm \tT(t)/s_e) S_{A,J},$$
where $x_{A,J,\GG}(t) = (2t)^{e_J}q(t)^{e_\GG-e_J}$ and
$q(t)$, $e(t)$, $\tT(t)$, $s_e$ are as defined above.

Note that $x_{\phi_0,\ov{e},e}(t)=q(t)$, so the good event pertaining to $Q(i)$ 
is $Q(i) = (1 \pm e(t)/s_e)(q(t) \pm \tT(t)/s_e)n^2$.
We also write $c(t) = aut(H)^{-1} \sum_T x_{\phi_T,J_T,\GG_T}(t)$,
where as above the sum is over quadruples $T=(a,b,c,d)$ such that
$ab$ and $cd$ are distinct edges of $H$.

Now we give an informal derivation of the differential equations satisfied by the functions $x_{A,J,\GG}(t)$,
which describe the main terms for the behaviour of the variables $X_{\phi,J,\GG}$.
We stress that this discussion does not constitute a proof of Theorem~\ref{track}; 
rather, it motivates the functions \( x_{A,J,\GG}(t) \) defined above,
and presages the proper proof given below,
in which the calculations we make here will play a central role.
For the sake of the discussion we ignore the error terms described by $e(t)$ and $s_e$,
and use the approximations $X_{\phi,J,\GG}(i) \approx x_{A,J,\GG}(t) S_{A,J}$,
so $Q(i) \approx q(t)n^2$ and $C_{uv}(i) \approx c(t) p^{e_H-2}n^{v_H-2} = c(t)p^{-1}$.
The system of differential equations will follow from the
approximation $x_{A,J,\GG}(t+s^{-1}) \approx x_{A,J,\GG}(t) + s^{-1}x'_{A,J,\GG}(t)$ 
and replacing changes $X_{\phi,J,\GG}(i+1)-X_{\phi,J,\GG}(i)$ by their expected value given $\mc{G}_i$.
Intuitively, although the change in a single step may be far from its expected value,
over many steps a `law of large numbers' will apply to the accumulated changes.
We also ignore two `pathological' behaviours that will need to be dealt with in Section \ref{sec:track}.
As an illustrative case, we start by counting open edges $|O(i)|=Q(i)/2$.
When we choose the edge $e_{i+1}$ we have
$$Q(i+1) = Q(i) - 1 - C_{e_{i+1}}(i) \approx q(t)n^2 - c(t) p^{-1}.$$
Since
$$Q(i+1) \approx q(t+1/s)n^2 \approx (q(t) + s^{-1}q'(t))n^2 = q(t)n^2 + p^{-1}q'(t)$$
we have the equation $q'(t) = - c(t)$.

To derive the differential equation for the general extension variable \( x_{A,J,\GG}(t) \), we 
write $$X_{\phi,J,\GG}(i+1)-X_{\phi,J,\GG}(i)=Y^+_{\phi,J,\GG}(i)-Y^-_{\phi,J,\GG}(i),$$
where $Y^+_{\phi,J,\GG}(i) \ge 0$ is the number of functions $f:V_\GG \to [n]$
in $\Xi_{\phi,J,\GG}(i+1) \sm \Xi_{\phi,J,\GG}(i)$, and
$Y^-_{\phi,J,\GG}(i) \ge 0$ is the number of functions $f:V_\GG \to [n]$
in $\Xi_{\phi,J,\GG}(i) \sm \Xi_{\phi,J,\GG}(i+1)$.
The term $Y^+_{\phi,J,\GG}(i)$ has contributions corresponding to each edge $e$ of $J$.
A function $f$ in $\Xi_{\phi,J \sm e,\GG}(i)$ will be counted by $Y^+_{\phi,J,\GG}(i)$
if the process chooses the edge $e_{i+1}$ equal to $f(e)$.
Since $e_{i+1}$ is chosen uniformly at random among $Q(i)/2$ open edges,
we can estimate
$$\mb{E}(Y^+_{\phi,J,\GG}(i) \vert \mc{G}_i) \approx
 \frac{2}{Q(i)} \sum_{e \in J} X_{\phi,J \sm e,\GG}(i)
\approx \frac{ 2p^{-1}S_{A,J}}{q(t)n^2} \cdot \sum_{e \in J} x_{A,J \sm e,\GG}(t).$$
The term $Y^-_{\phi,J,\GG}(i)$ has contributions corresponding to each edge $e$ of $\GG \sm J$.
A function $f$ in $\Xi_{\phi,J,\GG}(i)$ will be counted by $Y^-_{\phi,J,\GG}(i)$
if the process either chooses the edge $e_{i+1}$ equal to $f(e)$ or $f(e)$ becomes closed,
i.e. $f(e) \in C(i+1)$.
Thinking of $e_{i+1}$ as an ordered pair, the number of choices is
$2 + C_{f(e)}(i)$, each occurring with probability $Q(i)^{-1}$.
Therefore
$$\mb{E}(Y^-_{\phi,J,\GG}(i) \vert \mc{G}_i)
= \frac{1}{Q(i)}\sum_{e \in \GG \sm J} \sum_{f \in \Xi_{\phi,J,\GG}(i)} (2+C_{f(e)}(i))
\approx (e_\GG - e_J) \frac{c(t)p^{-1} x_{A,J,\GG}(t){S_{A,J}  }}{q(t) n^2}.$$
On the other hand, we have
\begin{align*}
Y^+_{\phi,J,\GG}(i)-Y^-_{\phi,J,\GG}(i) & = X_{\phi,J,\GG}(i+1)-X_{\phi,J,\GG}(i)
\approx (x_{A,J,\GG}(t+s^{-1})-x_{A,J,\GG}(t))S_{A,J} \\
& \approx s^{-1} x'_{A,J,\GG}(t) S_{A,J}
\end{align*}
so we have the equation
\begin{equation}
\label{eq:diff}
q(t) x'_{A,J,\GG}(t) = 2\sum_{e \in J} x_{A,J \sm e,\GG}(t) - (e_\GG - e_J)c(t)x_{A,J,\GG}(t).
\end{equation}
Note that the equation \( q'(t) = c(t) \) derived above is simply a special case of (\ref{eq:diff}).

To solve these equations we use the substitution $x_{A,J,\GG}(t) = q(t)^{e_\GG-e_J} z_\ell(t)$,
where we will see that the functions $z_\ell(t)$ can be parameterised by a single number $\ell=e_J$.
Then, since $q'(t)=-c(t)$, we have
$q(t)x'_{A,J,\GG}(t) = q(t)^{e_\GG-e_J+1} z'_\ell(t) - c(t) (e_\GG - e_J) q(t)^{e_\GG-e_J} z_\ell(t)$,
which also equals
$$2\sum_{e \in J} x_{A,J \sm e,\GG}(t) - (e_\GG - e_J)c(t)x_{A,J,\GG}(t)
= 2\ell q(t)^{e_\GG-e_J+1}z_{\ell-1}(t) - (e_\GG - e_J)c(t) q(t)^{e_\GG-e_J} z_\ell(t).$$
We deduce that $z'_\ell(t) = 2\ell z_{\ell-1}(t)$.
Now we use the initial conditions that $x_{A,J,\GG}(0)$ is equal to $1$ if $e_J=0$,
otherwise $0$ (e.g. \( q(0)=1 \)).  So $z_0(0)=1$ and $z_\ell(0)=0$ for $\ell>0$.
We obtain the solution $z_\ell(t)=(2t)^\ell$.
Also $q'(t)=-c(t)=-aut(H)^{-1} \sum_T x_{\phi_T,J_T,\GG_T}(t)=-aut(H)^{-1} 4e_H(e_H-1) q(t)(2t)^{e_H-2}$.
Integrating and substituting we conclude that
\begin{align*}
q(t) & = e^{-2e_H aut(H)^{-1} (2t)^{e_H-1}} \\
x_{A,J,\GG}(t) & = (2t)^{e_J} e^{-2(e_\GG-e_J)e_H aut(H)^{-1} (2t)^{e_H-1}} = (2t)^{e_J} q(t)^{e_\GG-e_J}.
\end{align*}

\begin{quote}
\nib{Remark.} 
As discussed above, we expect these random variables to evolve as they
do in the unconstrained random graph \( G(n,i) \).  Thus it is natural to compare the 
process $G(t)$ at time $t$ to the random graph $G(n,\rho)$, where
$\rho n^2/2 = i = tpn^2$, i.e. $\rho = 2tp$.  In \( G(n,\rho) \) we can
define open/closed pairs and the variables $X_{\phi,J,\GG}(i)$. 
For any ordered pair $uv$ in $[n]$,
edge $ab$ of $H$ and function $f:V_H \to [n]$ with $f(a)=u$, $f(b)=v$ the edges
of $f(H \sm ab)$ will all be present in $G(n,\rho)$ with probability $\rho^{e_H-1}$.
(For the purpose of this discussion we ignore the negligible contributions 
from functions $f$ that are not injective.)
Given $uv$, there are $2e_H n^{v_H-2}$ such functions $f:V_H \to [n]$, corresponding
to $2e_H aut(H)^{-1} n^{v_H-2}$ distinct sets of edges. 
The probability that $uv$ is open should be approximately
$$(1-\rho^{e_H-1})^{2e_H aut(H)^{-1} n^{v_H-2}} \approx
\exp \left( -(2tp)^{e_H-1}2e_H aut(H)^{-1} n^{v_H-2} \right) = q(t).$$
Similar reasoning applies to general extension variables, and the equations we
derived above agree with the corresponding equations for \( G(n,\rho) \).  
(See Spencer \cite{Sp2} for results on extension variables in this model.)
We could use this correspondence as the starting point of our discussion and as a heuristic
for the trajectories our variables follow, but this would not provide any
insight into how to prove that our random variables actually
follow the given trajectories. As we noted above, the calculations in this section
play a central role in the proof of Theorem~\ref{track}.
\end{quote}

\section{Strictly balanced graphs and balanced extensions} \label{sec:bal}

In this section we obtain some basic properties of our fixed strictly $2$-balanced graph $H$.
We also introduce a more general concept of strictly balanced extensions, and discuss the
manner by which arbitrary extensions can be decomposed into a series of such extensions.
First we recall the relevant definitions. We suppose that $H$ is {\em strictly $2$-balanced}, 
in the sense that $v_H, e_H \ge 3$ and $\frac{e_H-1}{v_H-2} > \frac{e_K-1}{v_K-2}$
for all proper subgraphs $K$ of $H$ with $v_K \ge 3$.
We also fix the parameter $$p = n^{-\frac{v_H-2}{e_H-1}}.$$
For any graph $\GG$ we define the {\em scaling} of $\GG$ to be $S_\GG =n^{v_\GG}p^{e_\GG}$.
The condition that $H$ is strictly $2$-balanced can be also be written as
$S_K > S_H$ for all subgraphs $K$ of $H$ with $2 < v_K < v_H$,
since $S_H = n^{v_H}p^{e_H} = pn^2$ and
$$S_K/S_H = n^{v_K-2}p^{e_K-1} =
n^{(e_K-1)\left( \frac{v_K-2}{e_K-1} - \frac{v_H-2}{e_H-1} \right)} > 1.$$
Note that the scaling $S_\GG$ is always an integer power of $n^{1/(e_H-1)}$.
It follows that the inequality $S_\GG >1$ actually implies $S_\GG \ge n^{1/(e_H-1)}$
and similarly that $S_\GG<1$ implies $S_\GG \le n^{-1/(e_H-1)}$.

The following lemma collects some simple properties of $H$ and $p$.

\begin{lemma} \label{sb} $ $
\begin{itemize}
\item[(i)] If $d$ is the largest integer for which $np^{d-1}>1$ then $H$ has minimum degree at least $d$.
\item[(ii)] We have $p>1/n$, and so $H$ has minimum degree at least $2$.
\item[(iii)] $H$ is a $2$-connected graph, and if $\{x,y\}$ is a cutset then $xy \notin E_H$.
\end{itemize}
\end{lemma}

\nib{Proof.}
First note that $H$ cannot have a vertex $v$ of degree at most $d-1$:
otherwise $S_H/S_{H \sm v} = np^{d(v)} > 1$,
which contradicts the fact that $H$ is strictly $2$-balanced.
We deduce that $H$ has minimum degree at least $1$.
Next, suppose for a contradiction that $p \le 1/n$.
Then $e_H \le v_H-1$.
However, for every connected subgraph $K$ of $H$ we have $e_K \ge v_K-1$,
so $\frac{e_K-1}{v_K-2} \ge 1 \ge \frac{e_H-1}{v_H-2}$,
which contradicts the definition of $H$ being strictly $2$-balanced.
Therefore $p>1/n$.
Now suppose for a contradiction that $H$ is not $2$-connected.
Then we can write $V_H=X \cup Y$ so that $E_H = E_{H[X]} \cup E_{H[Y]}$
and $|X \cap Y|=1$. Then $S_{H[X]}S_{H[Y]}=nS_H$, so without loss of generality
we have $S_{H[X]} \le (nS_H)^{1/2}$, and since $S_H=pn^2$
we have $S_{H[X]}/S_H \le (n/S_H)^{1/2} = (1/pn)^{1/2} < 1$.
This contradicts $H$ being strictly $2$-balanced, so $H$ is $2$-connected.
Finally, suppose that $\{x,y\}$ is a cutset, but that $xy \in E_H$.
Write $V_H=X \cup Y$ so that $E_H = E_{H[X]} \cup E_{H[Y]}$
and $X \cap Y=\{x,y\}$. Then $S_{H[X]}S_{H[Y]}=pn^2S_H=(pn^2)^2$,
so without loss of generality $S_{H[X]} \le pn^2 = S_H$. But this
contradicts $H$ being strictly $2$-balanced, so $xy \notin E_H$.
\qed

Recall that if $\GG$ is a graph and $A \sub V_\GG$
we define the {\em scaling} of the pair $(A,\GG)$ to be
$$S_{A,\GG} = p^{e_\GG-e_{\GG[A]}} n^{v_\GG-|A|}.$$
Note that $S_{A,\GG} = S_\GG/S_{\GG[A]}$.
Also, for any $A \sub B \sub V_\GG$ we have
$S_{B,\GG} = S_\GG/S_{\GG[B]} = S_\GG/S_{\GG[A]} \cdot S_{\GG[A]}/S_{\GG[B]} = S_{A,\GG}/S_{A,\GG[B]}$.
We say that $(A,\GG)$ is {\em strictly balanced} if
for any $A \subn B \subn V_\GG$ we have
$S_{A,\GG} < S_{A,\GG[B]}$, or equivalently $S_{B,\GG} < 1$.
For example, we can again rephrase our assumption that $H$ is strictly $2$-balanced
to say that for any edge $e=ab$ of $H$, with $A=\{a,b\}$ the pair $(A,H)$ is strictly balanced.
Indeed, $S_{A,H} = p^{e_H-1}n^{v_H-2} = 1$,
and for $A \subn B \subn V_H$ we have
$S_{B,H} = S_H/S_{H[B]} < 1$.

We will apply results on strictly balanced extensions to arbitrary pairs $(A,\GG)$
using the {\em extension series} $A = B_0 \subn B_1 \subn \cdots \subn B_d = V_\GG$ 
of $(A,\GG)$, which we construct by the following rule. If $(B_i,\GG)$ is not strictly balanced
then $B_{i+1}$ is chosen to be a minimal set $C$ with $B_i \subn C \subn V_\GG$
that minimises $S_{B_i,\GG[C]}=n^{|C|-|B_i|}p^{e_{\GG[C]}-e_{\GG[B_i]}}$,
otherwise we choose $B_d = B_{i+1} = V_\GG$.
For more compact notation we also write $S^A_i(\GG)=S_{B_i,\GG[B_{i+1}]}$.
We note the following properties of extension series.
\begin{itemize}
\item
$(B_i,\GG[B_{i+1}])$ is strictly balanced.
\item
For $i \ge 1$ we have
$S^A_i(\GG) = S_{B_i,\GG[B_{i+1}]} = S_{B_{i-1},\GG[B_{i+1}]}/S_{B_{i-1},\GG[B_i]} \ge 1$.
Therefore the sequence $S_{A,\GG[B_i]} = \prod_{j=0}^{i-1} S^A_j(\GG)$ is non-decreasing.
However, it is not necessarily true that the sequence of successive factors $S^A_i(\GG)$
is non-decreasing. For example, consider the $K_7$-free process,
where $p=n^{-1/4}$, and let $\GG=K_4$. Choosing $A$ of size $2$
we have $\GG[B_0]=K_2$, $\GG[B_1]=K_3$, $\GG[B_2]=K_4$
with $S^A_0(\GG) = np^2 = n^{1/2}$ and $S^A_1(\GG) = np^3 = n^{1/4}$.
\item
It is possible that $S_{A,\GG}<1$ but some factors $S^A_i(\GG)$ are greater than $1$.
For example, consider the $C_5$-free process,
where $p=n^{-3/4}$, and let $\GG$ be the graph consisting of $K_4$
plus an isolated vertex. Choosing $A$ to be $2$ vertices of the $K_4$
we have $\GG[B_0]=K_2$, $\GG[B_1]=K_4$, $\GG[B_2]=\GG$,
so $S^A_0(\GG) = n^2p^5 = n^{-7/4}$, $S^A_1(\GG)=n$ and $S_{A,\GG}=n^{-3/4}$.
\end{itemize}

\section{Union bounds}

In this section we collect some useful properties of the $H$-free process,
assuming that the good events $\mc{G}_i$ hold.
Recall that on $\mc{G}_i$ we have $Q(i)=(1 \pm e(t)/s_e)(q(t) \pm \tT(t)/s_e)n^2$,
and $q(t) = \exp\left(-\Theta(t^{e_H-1})\right)$,
where the constant in the $\Theta$-notation depends only on $H$.
We analyse the process up to time $t_{\max} = m/s = \mu (\log n)^{1/(e_H-1)}$,
and choose $\mu>0$ sufficiently small so that $e(t), q(t)^{-V} < n^\eps$.
Since $s_e = n^{1/2e_H-\eps}$ we have $Q(i) > n^{2-\eps}$ (say) for $i \le m$.
The following lemmas use this lower bound for $Q(i)$ and union bound estimates.
We will state the bounds at time $m$, but they also hold at any time $i \le m$ by monotonicity.
Our first lemma bounds the probability that $G(m)$ contains some fixed graph $F$.

\begin{lemma} \label{first-moment}
For any fixed graph $F$ on $[n]$, the probability that $\mc{G}_m$ holds and
$G(m)$ contains $F$ is at most $p^{e_F} n^{2e_F \eps}$.
\end{lemma}

\nib{Proof.} We take a union bound over
all choices of steps $1 \le i_1, \cdots, i_{e_F} \le m$ where the $j$th edge of $F$ is
chosen as the edge $e_{i_j}$ added to form $G(i_j)$ from $G(i_j-1)$. Since edges are
chosen uniformly at random from at least $n^{2-\eps}$ options, each choice has
probability at most $n^{-(2-\eps)}$ conditional on the history of the process.
Therefore $\mb{P}(F \sub G(m)) \le m^{e_F} n^{-(2-\eps)e_F} < p^{e_F} n^{2e_F \eps}$, say,
since $m = \mu (\log n)^{1/(e_H-1)} pn^2$. \qed

Given sets $A,B \sub [n]$, write $e(A,B)$ for the number of edges in $G(m)$
that have one endpoint in $A$ and the other in $B$. Our next lemma gives a bound for
$e(A,B)$ holding with high probability for all choices of $A,B$ of specified size.

\begin{lemma} \label{edges}
For any $a,b \ge 1$, the probability $p_{a,b}$
that $\mc{G}_m$ holds and there exist sets $A, B \sub [n]$
such that $|A|=a$, $|B|=b$ and $e(A,B) \ge \max\{4\eps^{-1}(a+b),pabn^{2\eps}\}$
satisfies $p_{a,b} < n^{-(a+b)}$.
\end{lemma}

\nib{Proof.} Write $x = \max\{ 4 \eps^{-1}(a+b),pabn^{2\eps}\}$.
We take a union bound over $\binom{n}{a}$ choices for $A$,
$\binom{n}{b}$ choices for $B$, at most $\binom{ab}{x}$ ways to choose
$x$ pairs with one endpoint in $A$ and the other in $B$,
and less than $m^x$ choices of steps $1 \le i_1 < \cdots < i_x \le m$
in which to choose these pairs as edges of the process.
Since edges are chosen uniformly at random from at least $n^{2-\eps}$ options, each choice has
probability at most $n^{-(2-\eps)}$ conditional on the history of the process.
Therefore we can estimate the probability by
$p_{a,b} < \binom{n}{a}\binom{n}{b}\binom{ab}{x}m^xn^{-(2-\eps)x}$.
Since $m = \mu (\log n)^{1/(e_H-1)} pn^2$, we have
\begin{eqnarray*}
\log p_{a,b}
& < & a(\log(n/a) + 1) + b(\log(n/b) + 1) \\
&   & + \ x(\log(ab/x)+1+\log (pn^\eps) + \log \mu + (e_H-1)^{-1}\log\log n) \\
& < & (a + b - \eps x/2)\log n,
\end{eqnarray*}
since $x \ge pabn^{2\eps}$ and $n$ large imply that
$-(\log(ab/x) + \log (pn^\eps)) \ge \eps \log n \gg \log\log n$.
Since $x \ge 4\eps^{-1}(a+b)$ the stated bound follows. \qed

For $A \sub [n]$ let $D_{A,d}$ be the set of vertices $v$ such
that $|N_{G(m)}(v) \cap A| \ge d$, i.e. in $G(m)$,
$v$ has at least $d$ neighbours in $A$. We conclude this section
by applying the previous lemma to give an upper bound for $D_{A,d}$.

\begin{lemma} \label{degree}
For any $8\eps^{-1} \le d \le a \le dp^{-1}n^{-2\eps}$,
the probability that $\mc{G}_m$ holds and there exists $A \sub [n]$ with $|A| = a$
and $|D_{A,d}| \ge 8\eps^{-1}d^{-1}a$ is at most $n^{-a}$.
\end{lemma}

\nib{Proof.}
Set $B=D_{A,d}$, $b=|B|$ and consider the event that $b \ge 8\eps^{-1}d^{-1}a$.
Since $e(A,B) \ge db$ and $d \ge 8\eps^{-1}$ we have
$e(A,B) - 4\eps^{-1}b \ge db/2 \ge 4\eps^{-1}a$.
Also, the bound $a \le dp^{-1}n^{-2\eps}$ implies
that $e(A,B) \ge db \ge pabn^{2\eps}$.
By Lemma \ref{edges} this event has probability at most
$n^{-(a+b)} \le n^{-a}$. \qed

\section{Counting extensions}

In this section we see how to obtain general upper bounds on
extension variables, assuming that the good events $\mc{G}_i$ hold.
We will state the bounds at time $m$, but they also hold at any time $i \le m$ by monotonicity.
Let $N_{\phi,J} = X_{\phi,J,J}(m)$: the number of extensions of a fixed embedding $\phi:A \to [n]$
to an embedding $f:J \to G(m)$, where $A \sub V_J$ is independent. Note that this is
an upper bound for $X_{\phi,J,\GG}(m)$. The following lemma gives a good estimate on
$N_{\phi,J}$ when the extension is strictly balanced.

\begin{lemma} \label{extend}
Suppose $(A,J)$ is strictly balanced and $\phi:A \to [n]$ is an injective map.
Let $\omega(n)$ be any function such that $\omega(n) \to \infty$ as $n \to \infty$.
On $\mc{G}_m$, with high probability we have
$N_{\phi,J} < S_{A,J} n^{4e_J \eps}$ if $S_{A,J} \ge 1$
and $N_{\phi,J} < \omega(n)$ if $S_{A,J} < 1$.
\end{lemma}

\nib{Proof.}
We start by estimating the maximum number of vertex-disjoint extensions of $\phi$ to an embedding of $J$.
Let $N'_{\phi,J}$ be the maximum number $s$ such that
there are embeddings $f_1,\cdots,f_s$ of $J$ in $G(m)$, all restricting to $\phi$ on $A$,
with $f_i(V_J \sm A)$ and $f_j(V_J \sm A)$ disjoint for all $1 \le i<j \le s$.
We can estimate $\mb{P}(N'_{\phi,J} \ge s)$ by a union bound over at most $s!^{-1} (n^{v_J-|A|})^s$
possible functions $f_1,\cdots,f_s$, where for each choice of functions, we can apply
Lemma~\ref{first-moment} to obtain an upper bound $p^{se_J} n^{2se_J \eps}$
on the probability that the graph $F = \cup_{i=1}^s f_i(J)$ is a subgraph of $G(m)$.
Therefore
$$\mb{P}(N'_{\phi,J} \ge s) \le s!^{-1} (n^{v_J-|A|})^s p^{se_J} n^{2se_J \eps}
< (3s^{-1}S_{A,J} n^{2e_J \eps})^s.$$
If $S_{A,J} \ge 1$ then we can set $s=S_{A,J} n^{3e_J \eps}$ to get a bound holding
with failure probability much less than $\exp\left(-n^\eps\right)$.
On the other hand, if $S_{A,J} =  p^{e_J} n^{v_J-|A|}<1$ then, since $p = n^{-\frac{v_H-2}{e_H-1}}$,
we in fact have $S_{A,J} \le n^{-1/(e_H-1)}$. Assuming that $\eps< (2e_Je_H)^{-1}$ we then
have $S_{A,J} n^{2e_J \eps}<1$, and we can set $s=\omega'(n)$
for any function $\omega'(n) \to \infty$ as $n \to \infty$ to get a bound holding
with failure probability much less than $n^{-C}$ for any constant $C>0$.

Now we argue by induction on $v_J-|A|$ to show the following bounds on $N_{\phi,J}$:
if $S_{A,J} \ge 1$ then $N_{\phi,J} < S_{A,J} n^{3e_J \eps} \omega'(n)^{2(v_J-|A|)}$
and
if $S_{A,J} < 1$ then $N_{\phi,J} < \omega'(n)^{2(v_J-|A|)}$.
Then we can choose $\omega'(n)^{2(v_J-|A|)} < \omega(n) < n^\eps$
to obtain the bounds required for the theorem.
Our base case is $v_J-|A|=1$, when we have $N_{\phi,J}=N'_{\phi,J}$,
and we can apply the bounds just shown for $N'_{\phi,J}$.

Next suppose $v_J-|A|>1$. We claim that for any embedding $f$ counted by $N_{\phi,J}$
there are at most $\omega'(n)^{2(v_J-|A|)-1}$ embeddings $f'$ counted by $N_{\phi,J}$ with
$f'(V_J \sm A) \cap f(V_J \sm A) \ne \emptyset$. To see this, consider any such $f'$
and let $B = \{b \in V_J: f'(b) \in f(V_J)\}$, so that $A \subn B \subn V_J$.
Let $\phi'$ be the restriction of $f'$ to $B$
and let $J' = J \sm E_{J[B]}$ be the graph obtained from $J$ by deleting all edges inside $B$.
Then, as noted above, $S_{B,J'} = S_{A,J}/S_{A,J[B]}$,
and since $(A,J)$ is strictly balanced we have $S_{B,J'}<1$.
By induction hypothesis we have $N_{\phi',J} < \omega'(n)^{2(v_J-|B|)}$.
Also, there are at most $v_J^{|B|-|A|} < v_J^{v_J}$ choices for $\phi'$,
so at most $v_J^{v_J} \omega'(n)^{2(v_J-|B|)}$ embeddings $f'$
corresponding to this set $B$.
Summing over all $A \subn B \subn V_J$ we obtain
at most $\omega'(n)^{2(v_J-|A|)-1}$ (say) such embeddings $f'$.

Finally, we can estimate $N_{\phi,J}$ by means of a maximum collection
$F = \{f_1,\cdots,f_s\}$ of vertex-disjoint extensions of $\phi$
(so $|F|=N'_{\phi,J}$).
Any extension $f$ counted by $N_{\phi,J}$ has a common image with some $f_i \in F$ outside of $A$,
and for each $f_i \in F$ we have at most $\omega'(n)^{2(v_J-|A|)-1}$ such embeddings $f$.
Therefore $N_{\phi,J} \le N'_{\phi,J} \omega'(n)^{2(v_J-|A|)-1}$.
If $S_{A,J} \ge 1$ then $N'_{\phi,J} < S_{A,J} n^{3e_J \eps}$
and so $N_{\phi,J} < S_{A,J} n^{3e_J \eps} \omega'(n)^{2(v_J-|A|)}$.
On the other hand, if $S_{A,J}<1$ then $N'_{\phi,J} < \omega'(n)$
and so $N_{\phi,J} < \omega'(n)^{2(v_J-|A|)}$.
This completes the proof. \qed

For general extensions $N_{\phi,J}$ may be considerably larger than $S_{A,J}$,
but the following lemma gives a useful bound.

\begin{lemma} \label{extend-general}
On $\mc{G}_m$, with high probability we have
$N_{\phi,J} < n^{4e_J \eps} \max_{A \sub B \sub V_J} S_{B,J}$.
\end{lemma}

\nib{Proof.}
Consider the extension series
$A = B_0 \subn B_1 \subn \cdots \subn B_d = V_J$.
We repeatedly apply Lemma \ref{extend} to bound the number of extensions
in each step of the series. At the first step we either have
$S^A_0(J)<1$ and so $N_{\phi,J[B_1]} < \omega(n)$
or $S^A_0(J) \ge 1$  and so $N_{\phi,J[B_1]} < S^A_0(J) n^{4e_{J[B_1]} \eps}$.
At subsequent steps $i \ge 1$ we have $S^A_i(J) \ge 1$,
so for each injection $\phi':B_i \to [n]$ we have
$N_{\phi',J_i[B_{i+1}]} < S^A_i(J) n^{4(e_{J[B_{i+1}]}-e_{J[B_i]}) \eps}$.
Multiplying these bounds and using $S_{A,J} = \prod_{i=0}^{d-1} S^A_i(J)$
gives a bound equal to either $n^{4e_J \eps} S_{A,J}$ when $S^A_0(J) \ge 1$
or $\omega(n) n^{4(e_J-e_{J[B_1]})\eps}S_{B_1,J}$ when $S^A_0(J)<1$.
By definition of the extension series, $\max_{A \sub B \sub V_J} S_{B,J}$
is either $S_{A,J}$ when $S^A_0(J) \ge 1$ or $S_{B_1,J}$ when $S^A_0(J)<1$.
Also, we may assume that $e_{J[B_1]} \ge 1$ (otherwise \(E_J\) is empty), so we can choose $\omega(n) < n^\eps$
to obtain the required bound. \qed

\begin{quote}
{\bf Remark.} In both of the preceding lemmas we can choose $\omega(n) = n^{c\eps}$
for some constant $c>0$ to make the failure probability exponentially small.
\end{quote}

We say that the pair $(A,J)$ is {\em dense} if $S^A_0(J) = S_{A,J[B_1]}\ge 1$
and {\em strictly dense} if  $S^A_0(J)  > 1$ (and so $S^A_0(J) \ge n^{1/(e_H-1)}$).
Since $S^A_i(J) \ge 1$ for $i \ge 1$,
for a dense pair we have $\max_{A \sub B \sub V_J} S_{B,J} = S_{A,J}$,
so the previous lemma gives an approximate upper bound of $S_{A,J}$ for $N_{\phi,J}$.
Note that if $(A,J)$ is strictly dense then so is $(A,J')$ for any subgraph $J'$ of $J$,
since we have $S_{A,J'[B]} \ge S_{A,J[B]} > 1$ for any $B$ with $A \subn B \sub V_J$.
The same argument shows that if $J$ is a subgraph of $H$ with $e_J \le e_H-2$
and $A = \{u,v\}$, where $uv \in E_H \sm E_J$, then $(A,J)$ is strictly dense.

We conclude this section by showing that adding an edge to a strictly dense pair
gives a significant improvement on the bound for $N_{\phi,J}$.

\begin{lemma} \label{add-edge}
Suppose that $(A,J)$ is a strictly dense pair, $a,b$ are vertices of $J$
with $ab \notin E_J$ and $\{a,b\} \not\sub A$, and $J' = J \cup \{ab\}$
is obtained by adding the edge $ab$ to $J$.
Then $\max_{A \sub B \sub V_{J'}} S_{B,J'} < S_{A,J}$,
and so on $\mc{G}_m$, with high probability we have $N_{\phi,J'} < n^{-1/(e_H-1)+4e_{J'}\eps} S_{A,J}$.
\end{lemma}

\nib{Proof.} Choose $B$ with $A \sub B \sub V_J$ maximising $S_{B,J'}$.
If $B=A$ we have $S_{B,J'} = pS_{A,J}$,
whereas if $B \ne A$ we have $S_{B,J'} \le S_{B,J} = S_{A,J}/S_{A,J[B]} < S_{A,J}$,
as $(A,J)$ is strictly dense. Either way we have
$S_{B,J'} \le n^{-1/(e_H-1)}S_{A,J}$, since it is an integer power of $n^{1/(e_H-1)}$,
so the bound on $N_{\phi,J'}$ follows from Lemma \ref{extend-general}. \qed

\section{Closure fidelity}

Recall that for an ordered pair $uv \in O(i)$,
we write $C_{uv}(i)$ for the set of ordered pairs $xy \in O(i)$ that would become closed,
i.e. belong to $C(i+1)$, if at time $i+1$ the process chooses $uv$ as the edge $e_{i+1}$.
By definition of $C(i+1)$ this means that adding $uv$ and $xy$ to $G(i)$ would create a copy of $H$.
Also, since $uv$ and $xy$ are open, any such copy of $H$ must use both $uv$ and $xy$.
In principle there could be many such copies of $H$, but we will show in this section that in fact
this is not the case, and moreover, by counting these copies of $H$ we obtain an
accurate estimate for the number of pairs closed by $uv$.

We frequently need to estimate the number of overlapping extensions of two 
pairs $(A_1,J_1)$ and $(A_2,J_2)$, so we will introduce some notation for this situation.
Recall that a graph $W$ is a {\em join} of two graphs $W_1$ and $W_2$ if it
has subgraphs $J_1$ isomorphic to $W_1$ and $J_2$ isomorphic to $W_2$
such that $V_W = V_{J_1} \cup V_{J_2}$ and $E_W = E_{J_1} \cup E_{J_2}$.
For convenient notation we use names for vertices in $J_1$ interchangeably
with their corresponding vertices in $W_1$, and similarly for $J_2$ and $W_2$.
Whenever we use this notation the sets $A_1$ and $A_2$ will be independent
and we will write $C = V_{J_1} \cap V_{J_2}$.

We need some further notation for describing the possibilities by which a pair $uv$ can
close a pair $xy$. There must be a subgraph $J$ obtained by deleting
two edges $ab$ and $cd$ from $H$ and an injective map
$f:V_H \to [n]$ such that $f(a)=u$, $f(b)=v$, $f(c)=x$, $f(d)=y$
and $f(e) \in E(i)$ for every edge of $J$.
The map $f$ is counted by $X_{\phi_T,J_T,\GG_T}(i)$, where
given such a quadruple $T=(a,b,c,d)$, we write $\GG_T = H \sm ab$, $J_T = H \sm \{ab,cd\}$
and define $\phi_T$ by $\phi_T(a)=u$ and $\phi_T(b)=v$.

For the sake of an argument needed in the proof of Lemma \ref{I-open}
we extend the definition of $C_{uv}(i)$ to allow the case when $uv \in C(i)$
is a closed pair: we define it as the number of pairs $xy$ such that adding
$uv$ and $xy$ to $G(i)$ creates a copy of $H$ containing both $uv$ and $xy$.

\begin{lemma} \label{closure}
With high probability, for every $1 \le i \le m$ and ordered pair $uv \in O(i) \cup C(i)$,
assuming $\mc{G}_i$,
we have $|C_{uv}(i)| = aut(H)^{-1} \sum_T X_{\phi_T,J_T,\GG_T}(i) \pm n^{-1/e_H} p^{-1}$,
where the sum is over quadruples $T=(a,b,c,d)$ such that
$ab$ and $cd$ are distinct (but not necessarily disjoint) edges of $H$.
\end{lemma}

\nib{Proof.}
Let $P$ be the set of ordered pairs $xy$ for which there exist (at least) two
embeddings $f_1,f_2$ of $H$ in $G(i) \cup \{uv,xy\}$ with $f_1(E_H) \ne f_2(E_H)$
such that both embedded copies $f_1(H)$ and $f_2(H)$ use the edges $uv$ and $xy$.
Given any $xy \in P$ we fix any two such embeddings $f_1$ and $f_2$.
Let $W$ be a graph isomorphic to $(f_1(H) \cup f_2(H)) \sm \{uv,xy\}$
and write $a,b,c,d$ for the vertices in $W$ corresponding to $u,v,x,y$ respectively.
Note that these are not necessarily distinct, but there are at least $3$ distinct vertices
in the list, since $\{u,v\} \ne \{x,y\}$. Let $\phi$ be the function defined by
$\phi(a)=u$ and $\phi(b)=v$. We bound $P$ by estimating, for all such $W$, the number
$N_{\phi,W}$ of embeddings of $W$ in $G(i)$ where $a$ is mapped to $u$ and $b$ to $v$.

There are two cases, according to whether or not we have $f_1(V_H) = f_2(V_H)$.
If $f_1(V_H) = f_2(V_H)$ then, since $f_1(E_H) \ne f_2(E_H)$,
$W$ is obtained from a subgraph $J = H \sm \{ab,cd\}$ of $H$ by adding at least one edge.
As noted above, $(ab,J)$ is strictly dense, and so by Lemma \ref{add-edge}
we have $N_{\phi,W} < n^{-1/(e_H-1)+4e_W\eps} p^{-1}$.
Now suppose that $f_1(V_H) \ne f_2(V_H)$.
We need to estimate $N_{\phi,W}$ where $W$ is the join of
$J_1 = f_1(H) \sm \{uv,xy\}$ and $J_2 = f_2(H) \sm \{uv,xy\}$.
With the above notation we have $A_1=A_2=\{a,b\}$,
and $C = V_{J_1} \cap V_{J_2}$ contains $\{a,b\}$ and $\{c,d\}$,
so $C \sm A_1$ and $C \sm A_2$ are non-empty.
Choose $B$ with $A_1 \cup A_2 \sub B \sub V_W$ maximising $S_{B,W}$
and write $B_1 = B \cap V_{J_1}$, $B_2 = B \cap V_{J_2}$.
We consider three subcases according to $B_1$ and $B_2$.
The first subcase is $B_1 \cup C \ne V_{J_1}$.
Then we have $S_{B_1 \cup C,J_1} = S_{B_1 \cup C,H} < 1$,
as $\{c,d\} \sub C$ and $H$ is strictly $2$-balanced.
Also $S_{B_2,J_2} \le S_{A_2,J_2}$, since $(A_2,J_2)$ is (strictly) dense,
so $S_{B,W} \le S_{B_2,J_2} S_{B_1 \cup C,J_1} < S_{A_2,J_2} = p^{-1}$.
The second subcase is $B_2 \cup C \ne V_{J_2}$,
when a similar argument gives
$S_{B,W} = S_{B_1,J_1} S_{B_2 \cup C,J_2} < S_{A_1,J_1} = p^{-1}$.
Finally, the third subcase is $B_1 \cup C = V_{J_1}$ and $B_2 \cup C = V_{J_2}$.
Then $V_{J_1} \sm (A_1 \cup C)$ and $V_{J_2} \sm (A_2 \cup C)$ are non-empty,
since $f_1(V_H) \ne f_2(V_H)$. Since $(A_1,J_1)$ is strictly dense we have
$S_{B_1,J_1} < S_{A_1,J_1} = p^{-1}$,
so $S_{B,W} = S_{B_1,J_1} S_{B_2 \cup C,J_2} = S_{B_1,J_1} < p^{-1}$.
In all cases we have $S_{B,W} < p^{-1}$, so $S_{B,W} \le n^{-1/(e_H-1)}p^{-1}$,
since it is an integer power of $n^{-1/(e_H-1)}$. Now Lemma \ref{extend-general}
gives $N_{\phi,W} < n^{-1/(e_H-1)+4e_W\eps} p^{-1}$.
Summing over less than $|V_H|^{2|V_H|}$ (say) choices of $W$
we obtain a bound $|P| \le n^{-1/(e_H-1/2)}p^{-1}$, say.

To finish the proof we calculate the number of ordered pairs $xy \notin P$ counted by $C_{uv}(i)$.
For each such pair $xy$ there is a unique copy $H^c$ of $H$ in $G(i) \cup \{uv,xy\}$.
For each quadruple $T=(a,b,c,d)$ in $H$ such that there is an isomorphism $f:H \to H^c$
with $f(a)=u$, $f(b)=v$, $f(c)=x$, $f(d)=y$ we count $xy$ by $X_{\phi_T,J_T,\GG_T}(i)$.
Also, any other such quadruple $T'=(a',b',c',d')$ and isomorphism $f':H \to H^c$
with $f'(a')=u$, $f'(b')=v$, $f'(c')=x$, $f'(d')=y$ corresponds to the automorphism
$f^{-1}f'$ of $H$, and this is a one-to-one correspondence.
Therefore we can estimate the number of ordered pairs $xy \notin P$ that close $uv$ by
$aut(H)^{-1} \sum_T (X_{\phi_T,J_T,\GG_T}(i) \pm |P|)$.
Including the pairs in $P$, we can estimate $|C_{uv}(i)|$ by
$aut(H)^{-1} \sum_T X_{\phi_T,J_T,\GG_T}(i) \pm n^{-1/e_H} p^{-1}$, say.
This completes the proof. \qed

Note that the extension variables which appear in Lemma \ref{closure}
are trackable: they satisfy condition (b) in the definition, since $uv \notin E(i)$.
Substituting the formulae 
$X_{\phi_T,J_T,\GG_T}(i) = (1 \pm e(t)/s_e)((2t)^{e_H-2}q(t) \pm \tT(t)/s_e)p^{-1}$
and recalling that $s_e = n^{1/2e_H-\eps} \ll n^{1/e_H}$ 
we obtain the following estimate.

\begin{coro} \label{closed-estimate}
With high probability, for every $1 \le i \le m$ and ordered pair $uv \in O(i) \cup C(i)$,
assuming $\mc{G}_i$, we have
$$|C_{uv}(i)|  = (1 \pm 2e(t)/s_e)(a_H (2t)^{e_H-2}q(t) \pm \tT(t)/s_e) p^{-1},$$
$$\mbox{ where } \qquad a_H=4e_H(e_H-1)/aut(H).$$
\end{coro}

\section{Martingale estimates: the differential equations method}

Our main tool for establishing concentration of random variables will be the
following versions of the Azuma-Hoeffding inequality, Lemmas 6 and 7 from \cite{B}.
First we need some definitions. Suppose we have a sequence of random variables
$X_0,X_1,\cdots$ and a filtration $\mc{F}_0 \sub \mc{F}_1 \sub \cdots$
(which will always be the natural filtration given by the process).
We say that the sequence $X_0,X_1,\cdots$ is a {\em martingale}
if $\mb{E}(X_{i+1}|\mc{F}_i)=X_i$ for $i \ge 0$.
We say it is a {\em submartingale} if $\mb{E}(X_{i+1}|\mc{F}_i)\ge X_i$ for $i \ge 0$
or a {\em supermartingale} if $\mb{E}(X_{i+1}|\mc{F}_i)\le X_i$ for $i \ge 0$.
We say that a sequence of random variables $X_0,X_1,\cdots$
is {\em $(\eta,N)$-bounded}, for some $\eta,N>0$, if
$X_i - \eta \le X_{i+1} \le X_i+N$ for all $i \ge 0$.
In our application below
we consider sequences of random variables \( A_0, A_1, \dots \) where
the difference sequence $D_i = A_{i+1}-A_i$ satisfies
$0 \le D_i \le N$ and $\mb{E}D_i = (1 \pm e_i)d_i$
for some $d_i \le \eta/2$ and a small error term $0<e_i<1$.
We will define
$A_i^+ = \sum_{j<i} (D_j - (1-e_j)d_j)$,
and
$A_i^- = \sum_{j<i} (D_j - (1+e_j)d_j)$.
Then each of $A^\pm_i$ is $(\eta,N)$-bounded,
$A^+_i$ is a submartingale and $A^-_i$ is a supermartingale.
We refer to $A^\pm_i$ as a {\em martingale pair with parameters $(\eta,N)$}.

\begin{lemma} \label{martingale-below}
Suppose $\eta \le N/10$, $m \ge 1$, $a>0$ and $A_0,A_1,\cdots$ is an $(\eta,N)$-bounded submartingale.
Then $\mb{P}(A_m \le A_0-a) \le e^{-a^2/3\eta m N}$.
\end{lemma}

\begin{lemma} \label{martingale-above}
Suppose $\eta \le N/10$, $m \ge 1$, $0 < a \le \eta m/10$
and $A_0,A_1,\cdots$ is an $(\eta,N)$-bounded supermartingale.
Then $\mb{P}(A_m \ge A_0+a) \le e^{-a^2/3\eta m N}$.
\end{lemma}

We now come to the formulation of the differential equations method.
Although it is technically involved, the idea behind it is quite simple.
We have a collection of sequences of random variables, and would like
to prove that certain asymptotic approximations hold with high probability 
at each step of each sequence. The asymptotic formulae are heuristically
derived by considering the one-step expected changes in these variables.
We let $\mc{G}_i$ be the event all formulae hold up to step $i$.
If, conditional on $\mc{G}_i$, the expected change of a
random variable from step $i$ to step $i+1$ is close to what it should be for these formulae to hold,
and we also have a useful absolute bound for these one-step changes,
then we can apply martingale estimates to show that the event $\mc{G}_i$
indeed holds with high probability.
We recommend the survey of Wormald \cite{W} for an introduction to this method,
and a comparison of Lemma~\ref{de} below with Theorem 5.1 in Wormald \cite{W} may be helpful.
We also note that Seierstad \cite{Se1, Se2} has recently given improved large deviation bounds and
a central limit theorem for the method under certain general criteria.
One difference in our theorem is that we phrase our result in terms of
a known smooth solution to a system of differential equations, and thus
side-step the issue of the existence of a solution.
However, the important difference is in the hypothesis for the bounds
on the one-step changes of the variables: by using Lemmas \ref{martingale-below}
and \ref{martingale-above} we can make do with much weaker estimates
than those needed to apply the general result from \cite{W}.

\begin{quote}
{\bf Set-up for Lemma~\ref{de}.}
Suppose we have a stochastic graph process defined on the vertex set \( [n] \), where \(n\) is large.
Let \(r\) be a fixed positive integer, and for each \(j \in [r]\) let \( k_j, S_j \) 
be parameters (which can depend on \(n\)).
Suppose that 
for each \( j \in [r]\) and \( A \in \binom{[n]}{ k_j } \)
there is a sequence of random variables \( X_{j,A}(i) \), defined for \( i =0, \dots, m\) and 
measurable with respect to the underlying graph process.  We suppose further that
\[ X_{j,A}(i+1) - X_{j,A}(i) = Y^+_{j,A}(i) - Y^-_{j,A}(i), \]
where  $Y^+_{j,A}(i), Y^-_{j,A}(i) \ge 0$.  We relate these sequences of random variables to
functions on $[0,\infty)$ by introducing \( t = i/s \) for some function \( s = s(n) \) that goes to
infinity.  
We hope to find a collection \( x_j(t) \) of continuous functions such that
\[ X_{j,A}(i) \approx x_j(t) S_j \]
for all \( j \in [r], A \in \binom{[n]}{k_j} \) and \( i = 0, \dots, m \).
Note that in our application $i$ will be the number of edges that have been added, 
and we can think of \(s\) as the time-scaling for the underlying process.
We can think of $1 \le j \le r$ as the `type' of a random variable
and the set $A$ as giving its `position' in the graph. The parameter \( S_j \) is the 
size-scaling for the \(j\)-th type of random variable. 
\end{quote}

Now we will formally state our lemma. Note that for technical reasons we also allow
the introduction of an additional sequence $\mc{H}_i$ of high probability events.

\begin{lemma} \label{de}
Let  $0<\eps<1$ and $c,C > 0$ be constants, and suppose that for 
each \(j \in [r]\) we have a parameter $s_j = s_j(n)$,
and functions $x_j(t)$, $e_j(t)$, $\tT_j(t)$, $\gG_j(t)$
that are smooth and non-negative for $t \ge 0$. 
For $i^* = 1, \dots, m$ let $\mc{G}_{i^*}$
be the event that
$$X_{j,A}(i) = \left(1 \pm \frac{e_j(t)}{s_j} \right) \left(x_j(t) \pm \frac{\tT_j(t)}{s_j} \right) S_j$$ 
for all $1 \le i \le i^* $, $1 \le j \le r$ and $A \in \binom{[n]}{k_j}$. Suppose that also there is a decreasing 
sequence of events \( \mc{H}_i\), $1 \le i \le m$ such that \( \mb{P}( \mc{H}_m \mid \mc{G}_m ) \to 1 \) as $n \to \infty$,
and that the following conditions hold:
\begin{enumerate}
\item(trend hypothesis) When conditioning on \( \mc{G}_{i} \wedge \mc{H}_i \) we have
$$\mb{E}Y^{\pm}_{j,A}(i) = \left( y^{\pm}_j(t) \pm  \frac{h_j(t)}{4s_j} \right)  \frac{S_j}{s},$$
for all \( j \in [r]\) and \( A \in \binom{[n]}{k_j} \),
where \( y_j^{\pm}(t) \) and \( h_j(t) \) are smooth non-negative functions such that
$$x'_j(t) = y^+_j(t) - y^-_j(t) \ \ \ \ \ \ \ \
\text{ and } \ \ \ \ \ \ \ \ h_j(t)=(e_jx_j+\gG_j)'(t);$$
\item(boundedness hypothesis)
For each \( j \in [r] \), conditional on \( \mc{G}_i \wedge \mc{H}_i \) we have
\[  Y^{\pm}_{j,A}(i)  < \frac{ S_j}{ s_j^2 k_j n^\eps}; \]
\item(initial condition) for all \( j \in [r]\) we have \( e_j(0)=\gG_j(0)=0 \);
and \( X_{j,A}(0) = S_j x_j(0) \) for all \( A \in \binom{[n]}{ k_j} \);
\item We have $n^{3\eps} < s < m < n^2 $, $ s \ge 40C s_j^2 k_j n^\eps $,  $n^{2\eps} \le s_j < n^{-\eps}s$,
\begin{gather*}
\inf_{t \ge 0} \tT_j(t) + e_j(t)x_j(t)/2 - \gG_j(t)/2 > c, \\
\sup_{t \ge 0} |y^{\pm }_j(t)|<C, \ \ \ \ \ \ \ \  
\sup_{t \ge 0} |x'_j(t)|<C, \ \ \ \ \ \ \ \ \int_0^\infty |x''_j(t)| \ dt<C, \\
\sup_{t \ge 0} |h_j(t)|<n^\eps, \ \ \ \ \ \ \ \ \ \int_0^\infty |h'_j(t)| \ dt<n^\eps.
\end{gather*}
\end{enumerate}
Then \( \mb{P}( \mc{G}_m \wedge \mc{H}_m  ) \to 1 \) as $n \to \infty$.
\end{lemma}

\nib{Proof.}
On the event $\mc{G}_i \wedge \mc{H}_i$ we define
$$Y^{\pm_1 \pm_2}_{j,A}(i) = 
Y^{\pm_1}_{j,A}(i) - (y^{\pm_1}_j(t) \mp_2 h_j(t)/4s_j) S_j/s.$$
(Recall our convention that this is shorthand for $4$ separate sequences of variables,
one for each way of choosing signs for $\pm_1$ and for $\pm_2$.)
If any event $\mc{G}_i$ or $\mc{H}_i$ fails we define
all $Y^{\pm_1 \pm_2}_{j,A}(i')$ to be $0$ for $i'>i$. Define
$$Z^{\pm_1\pm_2}_{j,A}(i) = \sum_{i'=0}^{i-1} Y^{\pm_1\pm_2}_{j,A}(i'), \quad \quad N_j = \frac{ S_j}{ s_j^2 k_j n^\eps}
\quad \quad \mbox{ and } \quad \quad \eta_j = 4CS_j/s.$$
Using the bounds $|h_j(t)| < n^\eps$, $s_j > n^{2\eps}$, $|y^{\pm}_j(t)|<C$
we see that $Z^{+ \pm}_{j,A}(i)$ and $Z^{- \pm}_{j,A}(i)$
are martingale pairs with parameter $(\eta_j,N_j+\eta_j)$. For example
$Z^{++}_{j,A}(i+1)-Z^{++}_{j,A}(i)=Y^{++}_{j,A}(i)=Y^+_{j,A}-(y^+_j(t) - h_j(t)/4s_j)S_j/s$
is a submartingale by the trend hypothesis,
is bounded above by $N_j + n^{-\eps}S_j/4s < N_j + \eta_j$ by the boundedness hypothesis
and below by $-CS_j/s > -\eta_j$. (The other cases are similar.)

Next we need the Euler-Maclaurin summation formula (see \cite{Ap}), which is as follows.
Suppose $f(t)$ is a smooth function and $a$ is a natural number.
Then $I = \int_{0}^a f(i) di$ can be approximated by
$S = \frac{1}{2}f(0) + f(1) + \cdots + f(a-1) +  \frac{1}{2}f(a)$
with error $|S - I| < \int_0^a |f'(i)| \ di$.
We apply the formula to $f(i) = x'_j(t(i))$ for any $j  \in [r]$
and $a=i^*$ with $1 \le i^* \le m$. Write $t^*=i^*/s$.
Then
$$I = \int_0^{i^*} x'_j(t(i)) \ di = \int_0^{t^*} x'_j(\tau) s \ d\tau = s \left(x_j(t^*) - x_j(0) \right)$$
and
$$|S-I| < \frac{1}{s} \int_0^{i^*} |x''_j(t(i))| \ di = \int_0^{t^*} |x''_j(\tau)| \ d\tau < C,$$
so
$$ \left|x_j(t^*) - x_j(0) - \frac{1}{s} \sum_{i=0}^{i^*-1} x'_j(t(i)) \right|
< \frac{1}{s} \left(  \left| \frac{x'_j(0)}{2} \right|+ \left| \frac{ x'_j(t^*)}{2} \right| +\int_0^{t^*} |x''_j(\tau)| \ d\tau \right)
< \frac{3C}{s}.$$
We can rewrite this as
\begin{equation} \label{eq:x}
\frac{1}{s} \sum_{i=0}^{i^*-1} x'(t(i))S_j = \left(x_j(t^*)- x_j(0) \pm \frac{3C}{s} \right)S_j.
\end{equation}
Similarly, our assumptions on $h_j$ and the initial conditions \( e_j(0)=\gG_j(0)=0 \) give
$|e_j(t^*)x_j(t^*)+\gG_j(t^*) - \sum_{i=0}^{i^*-1} h_j(t(i))/s| < 3n^\eps/s$, 
which we can rewrite as
\begin{equation} \label{eq:h}
\sum_{i=0}^{i^*-1} h_j(t(i))/4s_j \cdot S_j/s = (e_j(t^*)x_j(t^*)+\gG_j(t^*) \pm 3n^\eps/s) S_j/4s_j.
\end{equation}

Now we will estimate the probability that any event $\mc{G}_i$ fails. We can restrict attention
to events where all $\mc{H}_i$ hold, as by assumption they all hold with high probability.
Fix $1 \le j \le k$, $A \in \binom{[n]}{k_j}$, $1 \le i^* \le m$, $t^*=i^*/s$. 
Consider the event that $i^*$ is the first step at which $ \mc{H}_{i^*}$ holds but $\mc{G}_{i^*}$ fails
and that it fails for the variable $X_{j,A}(i^*)$.
One possibility is that
$X_{j,A}(i^*) > (1 + e_j(t^*)/s_j)(x_j(t^*) + \tT_j(t^*)/s_j) S_j$.
By definition
\begin{align*}
X_{j,A}(i^*) - X_{j,A}(0) - \sum_{i=0}^{i^*-1} x'(t(i))S_j/s
& = \sum_{i=0}^{i^*-1} ( Y^+_{j,A}(i) - y^+_j(t)S_j/s - Y^-_{j,A}(i) + y^-_j(t)S_j/s )\\
& = Z^{+-}_{j,A}(i^*)-Z^{-+}_{j,A}(i^*) + 2 \sum_{i=0}^{i^*-1} h_j(t(i))/4s_j \cdot S_j/s.
\end{align*}
Applying equation (\ref{eq:x}) gives
$$Z^{+-}_{j,A}(i^*)-Z^{-+}_{j,A}(i^*) + 2 \sum_{i=0}^{i^*-1} h_j(t(i))/4s_j \cdot S_j/s
>  (e_j(t^*) x_j(t^*) + \tT_j(t^*) + \tT_j(t^*)e_j(t^*)/s_j - 3Cs_j/s) S_j/s_j.$$
Then equation (\ref{eq:h}), $n^{2\eps} < s_j < n^{-\eps}s$ 
and $\tT_j(t^*) + e_j(t^*)x_j(t^*)/2 - \gG_j(t^*)/2 > c$ give
$$Z^{+-}_{j,A}(i^*)-Z^{-+}_{j,A}(i^*) > (e_j(t^*) x_j(t^*)/2 - \gG_j(t^*)/2 + 
\tT_j(t^*) - (n^\eps+3Cs_j)/s)S_j/s_j > cS_j/2s_j.$$
We deduce that $Z^{+-}_{j,A}(i^*)> cS_j/4s_j$ or $Z^{-+}_{j,A}(i^*)<-cS_j/4s_j$.
Now we apply Lemmas \ref{martingale-below} and \ref{martingale-above} with $a = cS_j/4s_j$,
which is valid using our assumptions $s \ge 40Cs_j^2 k_j n^\eps  $, $ s_j > n^{2\eps}$ and $m>s$
which give $\eta_j < N_j/10$ and $a < \eta_j m/10$.
We deduce that these events have probability at most
$$\exp (- (cS_j/4s_j)^2 / 3\eta_j m (N_j+\eta_j))
< \exp(- 5k_j\log n) \ll \left|\binom{[n]}{k_j}\right|^{-1} n^{-3k_j},$$
say. 
A similar bound holds for the probability that
$X_{j,A}(i^*) < (1 - e_j(t^*)/s_j)(x_j(t^*) - c/s_j) S_j$,
when we have $Z^{--}_{j,A}(i^*)> cS_j/4s_j$
or $Z^{++}_{j,A}(i^*)<-cS_j/4s_j$.
Taking a union bound over $1 \le j \le r$, $A \in \binom{[n]}{k_j}$ and $1 \le i^* \le m$
completes the proof. \qed

\section{Trackable variables} \label{sec:track}

To apply Lemma \ref{de} to the extension variables $X_{\phi,J,\GG}(i)$,
we need to estimate the expected and maximum number of extensions that
may be created or destroyed in each step of the process.
In this section we establish a bound on the maximum number of extensions
created or destroyed; in other words, we verify the boundedness hypothesis.
Also, in anticipation of the expected change calculations needed for the trend hypothesis,
we show that two types of pathological subgraph configurations that 
could potentially spoil these calculations are suitably rare.  
More specifically, we show that, on the event \( \mc{G}_i \),
there are very few extensions in \( \Xi_{\phi,J,\GG} \) that contain a pair of open
pairs \( e,f \) such that the inclusion of one as an edge causes the other to become closed,
and very few extensions in \( \Xi_{\phi, J, \GG} \) for which there are two edges in \( \phi( E_\GG \sm E_J) \) 
that can both be closed by the addition of the same edge \( e_{i+1} \).
We stress that we obtain these bounds whenever the variable
is {\em trackable} (as defined in Subsection \ref{notate1}). 
In particular, this condition holds for the extension variables 
that track the open routes to \(H\) less an edge, the central 
variables in the proof of Theorem~\ref{track}. 

We begin with a technical lemma that amounts to showing that if \( X_{\phi, J , \GG} \) is trackable then there
are no `implicitly' closed edges in \( E_\GG \setminus E_J \).
\begin{lemma}
\label{lem:noH}
If \( X_{\phi,J,\GG}(i) \) is a trackable variable 
and \( uv \in E_\GG \setminus E_J \) then there does {\bf not} exist 
\( C \subseteq V_H \) with an injective embedding \( \psi: C \to V_\GG \) such that
\begin{enumerate}
\item $\psi(H[C])$ is a subgraph of the graph 
$\GG' = \GG \cup \left( \phi^{-1}(E(i)) \cap \binom{A}{2} \right)$ 
obtained from $\GG$ by adding the edges $ab$ for all $a,b \in A$ with $\phi(a)\phi(b) \in E(i)$,
\item for any vertex $v \in C$ with \( \psi(v) \not\in A \), every neighbour of $v$ in $H$
belongs to $C$, and
\item there is some edge $e$ in $H[C]$ with \( \psi(e) = uv \).
\end{enumerate}
\end{lemma}
\nib{Proof.}
Assume for a contradiction that $\psi$ is an embedding satisfying conditions (1-3) of the lemma.
Define $A' = \{v \in C: \psi(v) \in A\}$. We claim that $|A'| \ge 2$. 
This is clear if $H$ contains an edge $e$ with $\psi(e) \sub A$.
Otherwise, condition (1) implies that $C \ne V_H$, as $H$ is not a subgraph of $\GG$ by
definition of trackability. Then condition (2) implies that $A'$ disconnects $H$,
and since \(H\) is 2-connected we deduce that $|A'| \ge 2$.

Now let \(K \) be the graph obtained from \( H[C] \) by deleting all edges inside $A'$.
Now \(K\) is isomorphic to a subgraph of \( \GG\) by condition (1), so
$S_{A, \GG[A \cup \psi(C)]} \le S_{A',K}$. Also, $S_{A',K} = n^{|C|-|A'|} p^{e_H(A',C \sm A')}$
is equal to $S_{(V_H \sm C) \cup A',H}$ by condition (2). This in turn is at most $1$,
as $H$ is strictly balanced. We deduce that $S_{A, \GG[A \cup \psi(C)]} \le 1$.

Note also that $\psi(C)$ is not contained in $A$, as by condition (3) it
contains the edge $\psi(e)=uv$ of $\GG$. This rules out the possibility that
$(A,\GG)$ is strictly dense, so it remains to consider possibility (b) in the
definition of trackability. In this case we must have $S_{A, \GG[A \cup \psi(C)]} = 1$,
and so $S_{(V_H \sm C) \cup A',H}=1$, when the fact that $H$ is strictly balanced
implies that $C=V_H$, $|A'|=2$ and $A' \in E_H$. However, the existence of such an
embedding of $H$ in $\GG'$ is specifically ruled out by the definition of trackability,
so we have the required contradiction. \qed

Now we are ready to verify the boundedness hypothesis.
Following the notation of Lemma \ref{de} we write
$X_{\phi,J,\GG}(i+1)-X_{\phi,J,\GG}(i)=Y^+_{\phi,J,\GG}(i)-Y^-_{\phi,J,\GG}(i)$,
where $Y^+_{\phi,J,\GG}(i) \ge 0$
is the number of maps $f$ in $\Xi_{\phi,J,\GG}(i+1) \sm \Xi_{\phi,J,\GG}(i)$
and $Y^-_{\phi,J,\GG}(i) \ge 0$
is the number of maps $f$ in $\Xi_{\phi,J,\GG}(i) \sm \Xi_{\phi,J,\GG}(i+1)$.
Recall that $f:V_\GG \to [n]$ is counted by $X_{\phi,J,\GG}(i)$ if
$f(e) \in O(i)$ for every $e \in E_\GG \setminus E_J$,
$f(e) \in E(i)$ for every $e \in E_J$, and
$f$ restricts to $\phi$ on $A$.
Then $f$ will be counted by $Y^-_{\phi,J,\GG}(i)$ if there is at least one
$e \in E_\GG \sm E_J$ such that $f(e)$ either becomes closed at step $i+1$
or is the edge $e_{i+1}$ chosen by the process at step $i+1$.
Also, for each edge $e$ of $J$ and $f$ counted by $X_{\phi,J \sm e,\GG}(i)$,
$f$ {\em might} be counted by $Y^+_{\phi,J,\GG}(i)$ if $e_{i+1}=f(e)$.
(We will see below that $f$ may not actually be counted, but for the purpose
of an upper bound we do not need to take this into account here.)

\begin{lemma}[Boundedness hypothesis] \label{changes}
With high probability, for every $1 \le i \le m$,
assuming $\mc{G}_i$ and that $X_{\phi,J,\GG}(i)$ is trackable, we have
$Y^+_{\phi,J,\GG}(i) \le n^{-1/e_H} S_{A,J}$ and $Y^-_{\phi,J,\GG}(i) \le n^{-1/e_H} S_{A,J}$.
\end{lemma}

\nib{Proof.}
We start with the variable $Y^+_{\phi,J,\GG}(i)$.
Fix an edge $e=ab$ of $J$ and suppose the process chooses the edge $e_{i+1}=uv$ in step $i+1$.
Let $A' = A \cup \{a,b\}$, $J' = J \sm E_{J[A']}$ and
define $\phi':A' \to [n]$ agreeing with $\phi$ on $A$
and satisfying $\phi'(a)=u$, $\phi'(b)=v$.
Note that one of $a$ or $b$ may belong to $A$, but not both, as $A$ is independent in $J$.
Any $f$ counted by $Y^+_{\phi,J,\GG}(i)$ with $f(a)=u$ and $f(b)=v$
is counted by $X_{\phi',J',\GG}(i)$; we can bound this by $N_{\phi',J'}$,
which by Lemma \ref{extend-general} is at most
$N_{\phi',J'} < n^{4e_{J'} \eps} \max_{A' \sub B \sub V_{J'}} S_{B,J'}$.
Since $A \subn A'$ and $(A,J)$ is strictly dense we have
$\max_{A' \sub B \sub V_{J'}} S_{B,J'} \le n^{-1/(e_H-1)} S_{A,J}$.
Summing over all edges $e$ of $J$ we estimate
$Y^+_{\phi,J,\GG}(i) < n^{-1/e_H} S_{A,J}$.

Now consider the variable $Y^-_{\phi,J,\GG}(i)$.
Suppose the process chooses the edge $e_{i+1}=uv$ in step $i+1$.
Fix an edge $e$ of $\GG \sm J$.
We want to estimate the number of embeddings $f$ in $\Xi_{\phi,J,\GG}(i)$
for which $f(e)$ is either equal to $e_{i+1}$ or becomes closed in step $i+1$.
Since $(A,J)$ is strictly dense, Lemma \ref{add-edge} gives
an upper bound of $n^{-1/(e_H-1)+4(e_J+1)\eps}S_{A,J}$ on the number
of embeddings $f$ with $f(e)=e_{i+1}$. 

Next consider an embedding $f$
where $f(e)=xy$ becomes closed in step $i+1$.
Then there is an embedding $f_2$ of $H$ in $G(i) \cup \{uv,xy\}$.
Write $C' = f(V_J) \cap f_2(V_H)$ and identify the sets
$f^{-1}(C')$ and $f_2^{-1}(C')$ as a set $C$ on which $f$ and $f_2$ agree.
Then we have $f_2(a)=u$, $f_2(b)=v$ for some $a,b \in V_H$, and we have some
$c,d \in C$ with $f(c)=f_2(c)=x$, $f(d)=f_2(d)=y$, where $\{c,d\} \ne \{a,b\}$ 
and $\{c,d\} \not\sub A$ (since $A$ is independent in $\GG$).
Write $H' = H \sm \{ab,cd\}$ and let $W$ be the join of $J_1=J$ 
and $J_2=H'$ formed by identifying vertices in $C$
and removing any edges within $A' = A \cup \{a,b\}$.

We claim that $S_{B,W} \le n^{-1/(e_H-1)}S_{A,J}$ for all  $A' \sub B \sub V_W$.  Fix 
such a set \(B\) and
write $B_1 = B \cap V_{J_1}$ and $B_2=B \cap V_{J_2}$.
We have
\[ S_{B,W} = S_{B_1,J} \cdot S_{B_2 \cup C,H} \cdot p^\beta \]
where \( \beta \) is the number of edges in \( J_2 \) joining
\( B_2 \setminus C \) and \( C \setminus B_2 \).  Since \( (A,J) \) is strictly dense
we have \( S_{B_1,J} \le S_{A,J} \), with equality only if \(B_1 = A \).  Furthermore,
since \( \{a,b\} \cup C \) has at least 3 vertices, we have \( S_{ C\cup B_2,H} \le 1 \), 
with equality only if \( C \cup B_2 = V_H \). Thus we can restrict our attention to
the situation where \( B_1 = A \), \( B_2 \supset V_{J_2} \setminus V_{J_1} \) 
and \( \beta = 0 \). In this case we will use Lemma~\ref{lem:noH} to obtain
a contradiction. We view $C$ as a subset of $V_H$ and let $\psi$ be the identification
of $C$ with the subset of $V_\GG$ which is also called $C$. We can assume that
condition (1) is satisfied, as otherwise $f$ is an extension of $\phi$ to
an embedding of a supergraph of $J$ and then we have the required estimate on
$S_{B,W}$ by Lemma \ref{add-edge}. Also, $\beta=0$ gives condition (2),
and $f_2(cd)=xy=f(e)$ with $e \in E_\GG \sm E_J$ and $c,d \in C$,
which gives condition (3). Thus Lemma~\ref{lem:noH} shows that this case does not
actually arise. We deduce that $S_{B,W} \le n^{-1/(e_H-1)}S_{A,J}$.

Now applying Lemma \ref{extend-general} and
summing over all possibilities for $e$ and $W$ gives the required bound
$Y^-_{\phi,J,\GG}(i) < n^{-1/e_H} S_{A,J}$. \qed

Now we turn to two technical issues regarding the expected values of 
\( Y^+_{\phi, J, \GG}(i) \) and \( Y^-_{\phi, J, \GG}(i) \).  
We would like to approximate these using our estimates for extension variables.  
In the case of \( Y^+_{\phi, J, \GG}(i) \), our first approximation is that for each 
edge $e$ of $J$, an embedding $f$ counted by $X_{\phi,J \sm e,\GG}(i)$
should be counted by $Y^+_{\phi,J,\GG}(i)$ if $e_{i+1}=f(e)$.  
However, we need to account for the possibility that the addition of the edge \( e_{i+1} = f(e) \) closes
some edge \(f(e')\) where \( e' \in E_\GG \setminus E_J \).  
In the case of \( Y^-_{\phi, J, \GG} \), we sum \( C_{f(uv)}(i) \) over \( uv \in E_\GG \setminus E_J\)
to estimate the number of open edges \(xy\) such that choosing \( e_{i+1} = xy \) causes a given embedding \(f\) in 
\( \Xi_{\phi, J, \GG} \) to leave this set. However, we need to account for the possibility that
there could be edges \(uv,u'v' \in E_\GG \sm E_J \) such 
that \( C_{f(uv)}(i) \) and \( C_{f(u'v')}(i) \) have large intersection.  
We now establish two lemmas showing that these two `pathological' possibilities 
have a negligible impact.

\begin{lemma}[Creation fidelity]
\label{credo}
If \( X_{\phi,J,\GG} \) is a trackable variable then, with high probability 
on the event \( \mc{G}_i\), the number of 
extensions \( f \in \Xi_{\phi,J,\GG} \) with the property that there 
are distinct \( uv,xy \in E_\GG \setminus E_J \) such that \( G(i) \cup \{f(uv),f(xy)\} \) 
contains a copy of
\( H \) is at most \( n^{-1/e_H} S_{A,J} \).
\end{lemma}

\nib{Proof.}
Let  \(uv,xy \in E_\GG \setminus E_J \) be distinct and fixed.
Consider any graph \(W\) given by the join \( J \) and a copy of 
\( H \) less two edges, where \(uv\) and \(xy\) are identified with these missing edges.  
As in Lemma \ref{changes} it suffices to show that
$S_{B,W} \le n^{-1/(e_H-1)}S_{A,J}$ for all $A \sub B \sub V_W$.
The argument is almost identical to that in Lemma \ref{changes}.
With the same notation we again have
$S_{B,W} = S_{B_1,J} \cdot S_{B_2 \cup C,H} \cdot p^\beta$.
We again have \( S_{B_1,J} \le S_{A,J} \), with equality only if \(B_1 = A \).
Furthermore, in the current lemma we have \( u,v,x,y \in C \), so $|C| \ge 3$,
and \( S_{ C\cup B_2,H} \le 1 \), with equality only if \( C \cup B_2 = V_H \). 
Then Lemma~\ref{lem:noH} applies as before to complete the proof. \qed

\begin{lemma}[Destruction fidelity]
\label{destroy}
If \( uv,u'v' \in O(i) \) are distinct then, on \( \mc{G}_i \), we have 
\( |C_{uv}(i)  \cap C_{u'v'}(i)| \le  n^{-1/e_H} p^{-1} \) with high probability.
\end{lemma}
\nib{Proof.}
Let \(ab\) and \( cd \) be distinct edges of \( H \) and set \( H_1 = H \sm \{ab,cd\} \).  
Similarly, let \(a'b' \) and \(c'd'\) be distinct edges of \(H\) and set \( H_2 = H \sm \{a'b', c'd' \} \).  
Now let \(W\) be any join of \( H_1 \) and \(H_2\)
where \( c=c'\) and \( d=d'\) but \( ab \neq a'b'\).  
Set \(A = \{a,b\} \cup \{a',b'\}\). Then $|A| \ge 3$.
Appealing to Lemma~\ref{extend-general}, it
suffices to show \( S_{B,W}< p^{-1} \) for all \( A \subseteq B \subseteq V_W \).
Fix such a set \(B\). Similarly to before we have
\( S_{B,W} \le S_{B_1,H_1} S_{ C \cup B_2, H_2} p^{\beta_2} \),
where \( B_1 = B \cap V_{H_1} \), \( B_2 = B \cap V_{H_2} \), and 
\( C = V_{H_1} \cap V_{H_2} \) and \( \beta_2 \) is the number
of edges in $H_2$ joining \( B_2 \sm C \) and \( C \sm B_2 \).

Note that $c,d \in C$, so $S_{ C \cup B_2, H_2} = S_{C \cup B_2,H} \le 1$,
with equality only when $C \cup B_2 = V_H$.
Also, since $H_1$ is strictly dense we have $S_{B_1,H_1} \le 1/p$,
with equality only when $B_1 = \{a,b\}$.
Thus we obtain the desired inequality \( S_{B,W}< p^{-1} \),
except possibly in the case when $C \cup B_2 = V_H$,
$B_1 = \{a,b\}$ and $\beta_2=0$.
Also, the same argument reversing the roles of $H_1$ and $H_2$
shows that we obtain the desired inequality, 
except possibly in the case when $C \cup B_1 = V_H$,
$B_2 = \{a',b'\}$ and $\beta_1=0$, where \( \beta_1 \) is the number
of edges in $H_1$ joining \( B_1 \sm C \) and \( C \sm B_1 \).
Since $H$ is $2$-connected, the only remaining possibility is when
$V_{H_1}=V_{H_2}$. But then $S_{B,W} \le S_{A,H_1} < 1/p$,
as $H_1$ is strictly dense and $|A| \ge 3$.
Thus in all cases we have the desired inequality.
\qed

\section{Trajectory verification and Tur\'an bounds} \label{track-pf}

Now we use the above bounds and Lemma \ref{de} to prove Theorem \ref{track},
which shows that trackable extension variables are well described by
the differential equations given earlier in the paper.
It will then follow that the process does indeed continue until
at least time $t=t_{\max}=\mu(\log n)^{1/(e_H-1)}$,
i.e. $m = \mu (\log n)^{1/(e_H-1)} pn^2$ edges.
In particular, it will follow that variables counting
common neighbours of $d$-sets with $p^d n>1$ and variables counting
extensions from non-edge pairs to subgraphs of $H$ with at most $e_H-2$ edges
satisfy these equations. Then Corollary \ref{neighbours} is an immediate
consequence of the formulae for common neighbours. In particular,
when $d=1$ we deduce the minimum degree statement needed to prove Theorem \ref{turan}.
To prove Theorem \ref{turan} we will show that the good event $\mc{G}_m$ holds
with high probability, i.e. for every $i \le m$ and trackable extension variable
$X_{\phi,J,\GG}(i)$ corresponding to a triple in $\TT$, we have
$$X_{\phi,J,\GG}(i) = (1\pm e(t)/s_e)(x_{A,J,\GG}(t) \pm \tT(t)/s_e) S_{A,J},$$
where $x_{A,J,\GG}(t) = q(t)^{e_\GG-e_J}(2t)^{e_J}$
and $t, s_e, S_{A,J}, q(t), e(t), \tT(t)$ are as defined in Subsection \ref{notate1}.

\nib{Proof of Theorem \ref{track}.}
To apply Lemma \ref{de} we arbitrarily number the triples in $\TT$ by $1 \le j \le r$
and identify the extension variables $X_{\phi,J,\GG}(i)$
with the variables $X_{j,A}(i)$ appearing in the statement of the lemma.
We take $e_j(t)=e(t)$ and $\tT_j(t)=\tT(t)$ for all $1 \le j \le r$.
The event \( \mc{H}_i \) is the event that the estimates given in 
Lemmas~\ref{closure},~\ref{changes},~\ref{credo}~and~\ref{destroy} hold up to step $i$.
We will give values for the other parameters of the lemma later in this proof.

We start with the main step, which is checking the trend hypothesis.
For the expected one-step changes $\mb{E}[Y^\pm_{\phi,J,\GG}(i)|\mc{G}_i \wedge \mc{H}_i]$ we analyse
the error terms in our earlier heuristic derivation.
We start with the variable $Q(i)$, which counts the number of ordered pairs that
are open at step $i$.
Write $Q(i+1)-Q(i)=Q^+(i)-Q^-(i)$ with $Q^+(i), Q^-(i) \ge 0$.
Since $Q(i+1) = Q(i) - 1 - |C_{e_{i+1}}(i)|$ we have $Q^+(i)=0$ and $Q^-(i) = 1 + |C_{e_{i+1}}(i)|$.
Then Corollary \ref{closed-estimate} gives
$$Q^-(i) = 1+(1 \pm 2e(t)/s_e)(a_H (2t)^{e_H-2}q(t) \pm \tT(t)/s_e) p^{-1}.$$
We have $q'(t)=y^+_q(t)-y^-_q(t)$, where $y^+_q(t)=0$ for all $t$
and $y^-_q(t) = c(t) = a_H (2t)^{e_H-2}q(t)$. 
We also have $h_q(t)=(eq+\gG)'(t)$.
Now $e'(t) = P'(t)e^{P(t)}>W(t^{e_H-2}+1)e^{P(t)}$
and $q'(t)/q(t) = -a_H(2t)^{e_H-2}$, so since $W \gg V \gg e_H$ we have
$h_q(t)/y^-_q(t) > (V + Wa_H^{-1}(2t)^{-(e_H-2)})e^{P(t)}$ for $t>0$.
Since $s=pn^2$ and $\tT(t)<1$
we easily have the required condition for $Q^-(i)$, namely
$$Q^-(i) = (y^-_q(t) \pm h_q(t)/4s_e) n^2/s.$$
(We only need this estimate for $\mb{E}(Q^-(i)|\mc{G}_i \wedge \mc{H}_i)$, but actually
it always holds on the event $\mc{G}_i$.)

Now we check the trend hypothesis in the general case.
We write $X_{\phi,J,\GG}(i+1)-X_{\phi,J,\GG}(i)=Y^+_{\phi,J,\GG}(i)-Y^-_{\phi,J,\GG}(i)$.
The term $Y^+_{\phi,J,\GG}(i)$ has contributions corresponding to each edge $e$ of $J$.
A function $f$ in $\Xi_{\phi,J \sm e,\GG}(i)$ will be counted by $Y^+_{\phi,J,\GG}(i+1)$
if the process chooses the edge $e_{i+1}$ equal to $f(e)$ and this choice of \( e_{i+1} \) 
does not close any edge in $f(E_\GG \sm E_J)$.
Now $e_{i+1}$ is chosen uniformly at random among $Q(i)/2$ open edges, so appealing to Lemma~\ref{credo}
we can estimate
$$\mb{E}(Y^+_{\phi,J,\GG}(i) \vert \mc{G}_i \wedge \mc{H}_i)
= 2Q(i)^{-1} \sum_{e \in J} \left( X_{\phi,J \sm e,\GG}(i) \pm  n^{-1/e_H} S_{A,J \sm e} \right).$$
Now $X_{\phi,J \sm e,\GG}(i)= (1\pm e(t)/s_e)(x_{A,J \sm e,\GG}(t) \pm \tT(t)/s_e) S_{A,J \sm e}$.
Since $S_{A,J \sm e} = p^{-1}S_{A,J}$, $n^{-1/e_H} \ll 1/s_e$ and $\tT(t) \ge 1/2$ for $t \ge 0$
we estimate $\mb{E}(Y^+_{\phi,J,\GG}(i) \vert \mc{G}_i \wedge \mc{H}_i)$ as
$$2 ((1 \pm e(t)/s_e)(q(t) \pm \tT(t)/s_e)n^2)^{-1} \cdot e_J \cdot
(1\pm e(t)/s_e ) 
(q(t)^{e_\GG-e_J+1}(2t)^{{e_J}-1} \pm 2 \tT(t)/s_e) p^{-1}S_{A,J}.$$
We have $x'_{A,J,\GG}(t) = y^+_{A,J,\GG}(t) - y^-_{A,J,\GG}(t)$,
where $y^+_{A,J,\GG}(t) = 2{e_J} q(t)^{e_\GG-e_J}(2t)^{{e_J}-1}$
and $y^-_{A,J,\GG}(t) = a_H(e_\GG-e_J)q(t)^{e_\GG-e_J}(2t)^{{e_J}+e_H-2}$.
We also have $h_{A,J,\GG}(t)=(ex_{A,J,\GG}+\gG)'(t)$.
To establish the required bound, i.e.
$$\mb{E}(Y^+_{\phi,J,\GG}(i) \vert \mc{G}_i \wedge \mc{H}_i)
= (y^+_{A,J,\GG}(t) \pm h_{A,J,\GG}(t)/4s_e)S_{A,J}/s,$$
it suffices to show that
\begin{equation}
\label{eq:y+}
(1 \pm 4e(t)/s_e)(1 \pm 2 \tT(t)q(t)^{-1}/s_e)(1 \pm 2 \tT(t)(q(t)^{e_\GG-e_J+1}(2t)^{{e_J}-1})^{-1}/s_e)
\end{equation}
$$\subseteq \ 1 \pm (2{e_J} q(t)^{e_\GG-e_J}(2t)^{{e_J}-1})^{-1} h_{A,J,\GG}(t)/4s_e.$$
Setting \( x(t) = x_{A,J,\GG}(t) = (2t)^{e_J} q(t)^{e_\GG-e_J}\) we see that it is necessary to
establish that
\begin{equation}
\label{eq:poserr}
\frac{4 e_J e(t)x(t)}{t} + \frac{ 2 e_J \tT(t) x(t)}{t q(t)} + \frac{ 4 e_J \tT(t)}{q(t)} 
\end{equation}
is bounded above by
\begin{align*}
\frac{1}{4} \left( x(t) e'(t) \right. +  & \left.x'(t) e(t) + \gamma'(t) \right)\\
& = \frac{e^{P(t)} x(t)}{4} \left( W t^{e_H-2} + W  + 
\left(e_J/t  - (e_\GG- e_J)a_H (2t)^{e_H-2} \right) \frac{ e^{P(t)}-1}{ e^{P(t)}} \right) + \frac{ \gamma'(t)}{4} \\
& > \frac{e^{P(t)} x(t)}{4} \left( \frac{W}{2} t^{e_H-2} + W  \right) + \frac{ \gamma'(t)}{4}.
\end{align*}
Note that establishing this bound is in fact sufficient. To see this we observe that our
choice of \( \gamma(t) \) ensures that \( h_{A,J,\GG}(t) \) is bounded below by some constant 
(which is a function of \(V\)). Therefore the terms omitted in (\ref{eq:y+}) are $O(1/s_e)=o(1)$,
so do not cause the inequality to be violated when $n$ is sufficiently large.
Note also that we can assume that $e_J>0$, as otherwise \( Y^+_{\phi, J, \GG} = y^+_{A,J, \GG} = 0 \).
To verify the bound for \( t < 40V/W \) we note that \( x(t) \le 9t/4 \), as $e_J>0$, and therefore
(\ref{eq:poserr}) is at most \( 9 V e(t) + 15 V < 10 V e^{40 V} < \gamma'(t)/4 \).  
On the other hand, for \( t > 40V/W \) we note that the first two terms in (\ref{eq:poserr}) 
can each be bounded by \( W e^{P(t)} x(t)/10 \); the remaining term is bounded by \( \gG'(t)/4 > 5V \) for 
\( 40V/W < t < 1/(50V) \) and by \( e^{P(t)} x(t) \) for larger \(t\).

Next consider the term $Y^-_{\phi,J,\GG}(i)$, which has contributions 
corresponding to each edge $e$ of $\GG \sm J$.
A function $f$ in $\Xi_{\phi,J,\GG}(i)$ will be counted by $Y^-_{\phi,J,\GG}(i+1)$
if the process either chooses the edge $e_{i+1}$ equal to $f(e)$ or $f(e)$ becomes closed, i.e. $f(e) \in C(i+1)$.
Thinking of $e_{i+1}$ as an ordered pair, the number of choices is
$2 + |C_{f(e)}(i)|$, each occurring with probability $Q(i)^{-1}$.
Therefore, appealing to Lemma~\ref{destroy}, we have
$$\mb{E}(Y^-_{\phi,J,\GG}(i) \vert \mc{G}_i \wedge \mc{H}_i)
= Q(i)^{-1} \sum_{f \in \Xi_{\phi,J,\GG}(i)} \sum_{e \in \GG \sm J}
(2+|C_{f(e)}(i)|\pm n^{-1/e_H} p^{-1}).$$
We can estimate $|C_{f(e)}(i)|$ by Corollary \ref{closed-estimate},
so we estimate $\mb{E}(Y^-_{\phi,J,\GG}(i) \vert \mc{G}_i \wedge \mc{H}_i)$ as
\begin{align*}
& ((1 \pm e(t)/s_e)(q(t) \pm \tT(t)/s_e)n^2)^{-1} \cdot (e_\GG - e_J) \cdot
  (1\pm e(t)/s_e)(q(t)^{e_\GG-e_J}(2t)^{e_J} \pm \tT(t)/s_e) S_{A,J} \\
& \cdot
  (1 \pm e(t)/s_e \pm n^{-1/e_H} ) (a_H (2t)^{e_H-2}q(t)\pm 2 \tT(t)/s_e) p^{-1}.
\end{align*}
Now to establish the required bound, i.e.
$$\mb{E}(Y^-_{\phi,J,\GG}(i) \vert \mc{G}_i \wedge \mc{H}_i)
= (y^-_{A,J,\GG}(t) \pm h_{A,J,\GG}(t)/4s_e)S_{A,J}/s,$$
it suffices to show that
$$(1 \pm 4e(t)/s_e)(1 \pm 2\tT(t)q(t)^{-1}/s_e)(1 \pm \tT(t)(q(t)^{e_\GG-e_J}(2t)^{e_J})^{-1}/s_e)
(1 \pm 2 \tT(t)(a_H (2t)^{e_H-2}q(t))^{-1}/s_e) $$
$$\subseteq \ 1 \pm (a_H ( e_\GG - e_J) q(t)^{e_\GG-e_J}(2t)^{{e_J}+e_H-2})^{-1} h_{A,J,\GG}(t)/4s_e.$$
And this reduces to showing that
\begin{gather*}
4 a_H ( e_\GG - e_J)(2t)^{e_H-2} x(t)e(t) + \frac{  2 a_H ( e_\GG - e_J)(2t)^{e_H-2} x(t) \tT(t)}{q(t)}  \\ 
+\ a_H ( e_\GG - e_J)(2t)^{e_H-2} \tT(t) +  \frac{2( e_\GG - e_J)  x(t) \tT(t) }{ q(t)}
\end{gather*}
is bounded above by
$$  \frac{e^{P(t)} x(t)}{4} \left( \frac{W}{2} t^{e_H-2} + W  \right) + \frac{ \gamma'(t)}{4}.$$
This follows by estimates very similar to those given above for $Y^+_{\phi,J,\GG}(i)$.
We omit the details, except for remarking that is helpful to observe
that the term \( a_H( e_\GG - e_J) (2t)^{e_H-2} \) is bounded by \( \gamma'(t)/4 \) for \( t < 1/(50V) \).

This verifies the trend hypothesis of Lemma \ref{de}. To finish the proof we check the remaining conditions.
The boundedness hypothesis follows from
Lemma~\ref{changes} as we have \( n^{1/e_H} \gg n^{1/e_H - \eps} =  s_e^2 n^\eps \).  We have $|\TT|=r < V^{3V}$, $n^{2\eps} 
< s_e^2 n^\eps < n < pn^2 = s < m < n^2$ 
and $s_e = n^{1/2e_H - \eps} > n^{2\eps}$. The functions $x_{A,J,\GG}(t)$ and $y^\pm_{A,J,\GG}(t)$ all have
the form $F(t)e^{-Kt^{e_H-1}}$, where $F$ is a polynomial of degree at most $V+e_H$,
and $K$ and all coefficients in $F$ are non-negative and bounded above by $W$, say.
Here we can use $$\int_0^\infty t^a e^{-t} = a!  \qquad \mbox{and} \qquad
\sup_{t \ge 0} t^a e^{-t} = (a/e)^a \qquad \mbox{for} \ a \in \mb{N}$$
to see that $\sup_{t \ge 0} |y^{\pm}_{A,J,\GG}(t)|$, $\sup_{t \ge 0} |x'_{A,J,\GG}(t)|$ and $\int_0^\infty |x''_{A,J,\GG}(t)| \ dt$
are all bounded by some constant $C$ depending only on $W$.
Also, recall that $e(t)=e^{P(t)}-1$ with $P(t) = W(t^{e_H-1} + t)$, 
$h_{A,J,\GG}(t)=(ex_{A,J,\GG}+\gG)'(t)$, and $\gG(t)$ is a smooth increasing function 
such that $\gG(t)$ and $\gG'(t)$ are bounded by absolute constants.
The initial conditions $e(0)=\gG(0)=0$ hold. Since $t<t^* = \mu (\log n)^{1/(e_H-1)}$, 
by choosing $\mu$ sufficiently small we can ensure that $\sup_{t \ge 0} |h_{A,J,\GG}(t)|<n^\eps$
and $\int_0^\infty |h'_{A,J,\GG}(t)| \ dt<n^\eps$. 
Finally, we can choose $c=1/2$, since $\tT(t) = 1/2 + \gG(t)$,
so $\tT(t) + e(t)x(t)/2 - \gG(t)/2 > 1/2$.
\qed

\section{Counting small subgraphs} \label{sec:count}

In this short section we apply our results to count small subgraphs
in the $H$-free process and compare these counts to those known
for the $G(n,p)$ model. A rough summary is that the $H$-free process
looks very much like $G(n,p)$ from this perspective, except that
it does not contain any graphs that contain $H$. A more precise
description is given by Theorem \ref{count}, which we now prove.

\nib{Proof of Theorem \ref{count}.}
Statement (i) follows from Lemma~\ref{first-moment}, as \( \GG[B] \) does not 
appear in \( G(i) \) with high probability, and therefore 
\( \GG\) itself does not appear with high probability
(note that the failure probability here decays polynomially in $n$, not exponentially).
Statement (ii) follows from Theorem~\ref{track} applied to the
trackable variable $X_\GG(i)=X_{\emptyset, \GG, \GG}(i)$.
It remains to consider the case when $S_{\GG[B]} \ge 1$ for all $B \sub V_\GG$.
Form the extension series $\emptyset = B_0 \subn B_1 \subn \cdots \subn B_d = V_\GG$, 
as defined in Section~\ref{sec:bal}. We divide the \(m\) steps of the process into \(d\) equal intervals, 
and in the \(j\)th interval we show that with high probability there is an extension from 
a {\em fixed} copy of $\GG[B_{j-1}]$ (found in the previous interval) to a copy of $\GG[B_j]$.
By construction every step of the extension series is
strictly balanced, and our assumption in this case implies that the scalings
in each step satisfy $S_{B_{j-1},\GG[B_j]} \ge 1$.
Suppose that $\phi:B_{j-1} \to [n]$ is an embedding of $\GG[B_{j-1}]$ in $G((j-1)m/d)$.
If $S_{B_{j-1},\GG[B_j]} > 1$ then the variable $X_{\phi,\GG[B_j],\GG[B_j]}(i)$
is trackable, so the required extension exists by Theorem~\ref{track} 
(in fact there are many such extensions).
On the other hand, if $S_{B_{j-1},\GG[B_j]} = 1$ we can apply Theorem~\ref{track} 
to the trackable variables $X_{\phi,\GG[B_j] \sm e,\GG[B_j]}(i)$ 
with $e \in E_{\GG[B_j]} \sm E_{\GG[B_{j-1}]}$. Writing $a_j = e_{\GG[B_{j+1}]} - e_{\GG[B_j]}$
we can estimate the probability that in step $i$ the edge $e_i$ completes some
embedding of $\GG[B_j] \sm e$ for some $e$ to an embedding of $\GG[B_j]$
by $Q(i)^{-1} \sum_e X_{\phi,\GG[B_j] \sm e,\GG[B_j]}(i) \sim a_j(2t)^{ a_j -1}/(pn^2)$.
Since the length of each interval is $m/d \gg s = pn^2$ and $t \gg 1$
(ignoring the first half of the first interval, say)
we see that the required extension appears with high probability. \qed

\begin{quote} {\bf Remark.}
Our results for counting labelled copies of $\GG$ in the $H$-free process mirror
those obtained for the analogous counts in $G(n,p)$. However, rather more is known
in the $G(n,p)$ model, some of which is surveyed in Section VII of \cite{R}.
In the supercritical case Barbour, Karo\'nski and Ruci\'nski \cite{BKR}
gave a central limit theorem with estimates on the rate of convergence for the
appropriately normalised count. Spencer \cite{Sp2} analysed the critical case:
one of his results concerns the case when $\GG$ is strictly balanced,
when he obtains the asymptotic probability for $\GG$ to appear when $p$ is near the threshold.
It seems plausible that similar results may hold for the $H$-free process:
in the supercritical case one would need to extract distributional information 
from the differential equations method (along the lines of \cite{Se1}),
and in the critical case one would need a more accurate analysis of the above proof
(which seems to suggest a Poisson approximation). For the sake of 
brevity we do not pursue these possibilities here.
\end{quote}

\section{Smooth independence}

We have now shown that the $H$-free process continues until at least the time
$t_{\max} = \mu (\log n)^{1/(e_H-1)}$, when it has $m = \mu (\log n)^{1/(e_H-1)} pn^2$ edges.
In this section we describe an additional assumption (`smooth independence') on $H$,
under which we show that the independence number of the resulting graph is at most
$$\aA = 3\mu^{-1} (\log n)^{1-1/(e_H-1)} p^{-1}.$$
Since the independence number cannot increase when more edges are added,
we also have the same upper bound for the terminal graph of the process.
The main step of our proof will be to show that, for any set $I$ of size $\aA$,
with high probability we can track the number of open pairs contained within $I$:
at time $t$ there will be roughly $q(t)|I|^2$ open ordered pairs in $I$.
Then a simple union bound calculation will show that with high probability
$I$ is not independent at time $t_{\max}$.

To track the open pairs within a set $I$ we use Lemma \ref{de},
but we cannot simply apply the lemma directly, due to the possibility of closing
a large number of pairs in $I$ in a single step of the process.  Note that
in this application of Lemma~\ref{de} we will take \( k_j = \alpha \) 
and \( S_j = \alpha^2 \).  So we will not be able to achieve the boundedness
hypothesis in a useful way if we allow our process to close \( \alpha \) edges in
the set \(I\) in a single step (and this certainly is a possibility 
for many choices of \(H\)).
To deal with this, we say that the edge $e_i$ added in step $i$ is {\em $I$-good}
if it closes at most $n^{-5\eps}p^{-1}$ ordered pairs in $I$, otherwise $e_i$ is {\em $I$-bad}.
Then we say that a pair $uv$ in $I$ is {\em $I$-closed} at step $i$
if there is some step $i' \le i$ such that $e_{i'}$ is $I$-good 
and \( G(i') \cup \{uv\} \) contains a copy of \(H\).
If $uv$ in $I$ is not in $E(i)$ and not $I$-closed we say that it is {\em $I$-open} at step $i$.
Note that an $I$-closed pair is closed, 
but an $I$-open pair could be open or closed (but not an edge).
Let $Q_I(i)$ be the number of open ordered pairs in $I$ at step $i$
and $X_I(i)$ be the number of $I$-open ordered pairs in $I$ at step $i$.
We write $P_I \sub E(m)$ for the set of ordered edges at time $t_{\max}$ that are $I$-bad.
Then we say that $H$ has {\em smooth independence}
if with high probability $|P_I| < n^{-5\eps}p^{-1}$ for every set $I$ of size $\aA$.

Our first step is to apply Lemma \ref{de} to track the number
of $I$-open pairs in $I$.

\begin{lemma} \label{I-open}
If $H$ has smooth independence, then with high probability, for any set $I$ of size $\aA$,
the number of $I$-open ordered pairs in $I$ at step $i$
is $X_I(i) = (1 \pm e(t)n^{-2\eps})(q(t) \pm n^{-2\eps})\aA^2$.
\end{lemma}

\nib{Proof.}
We apply Lemma \ref{de} with $r=1$, $k_1=\aA$,
$X_{1,I}(i)=X_I(i)$ for $I \in \binom{[n]}{\alpha} $, $x_1(t) = q(t)$, 
$e_1(t)=e(t)$, $\gG_1(t)=\gG(t)$, $\tT_1(t)=\tT(t)$,
$s_1=n^{2\eps}$ and $S_1=\aA^2$.  We let \( \mc{H}_i \) be the event
that the estimates given by Theorem~\ref{track} hold up to step $i$
and that $|P_I| < n^{-5\eps}p^{-1}$ for every set $I$ of size $\aA$.

The main step is verifying the trend hypothesis of Lemma \ref{de}.
Note that adding an edge cannot create any new $I$-open pairs,
so we always have $Y^+_{1,I}(i)=0$. Now we calculate the expected one-step
change $\mb{E}(Y^-_{1,I}|\mc{G}_i \wedge \mc{H}_i)$. Recall that a pair $e$ becomes
closed at step $i+1$ if the process chooses the edge $e_{i+1}$ in $C_e(i)$
so a pair $e$ in $I$ becomes $I$-closed if is $I$-open and $e_{i+1}$
is chosen in $C_e(i) \sm P_I$. Also, if $e$ in $I$ is open as well as $I$-open
it may become an edge if the process chooses $e_{i+1}=e$.
Now $e_{i+1}$ is chosen uniformly among
$Q(i)$ open ordered pairs at step $i$, so
$$\mb{E}(Y^-_{1,I}|\mc{G}_i \wedge \mc{H}_i) = Q(i)^{-1} \sum_{e \in X_I(i)} (|C_e(i) \sm P_I| \pm 1).$$
(Here we also wrote $X_I(i)$ for the set of $I$-open pairs in $I$.)
Temporarily ignoring the error terms, this suggests the equation
$x'_1(t) = - q(t)^{-1}x_1(t)c(t)$, which has $q(t)$ as a solution,
explaining our choice of $x_1(t)$ above.
To account for the error terms, we estimate $Q(i)$ by Theorem \ref{track},
$C_e(i)$ by Corollary \ref{closed-estimate}, $X_I(i)$ by the fact that we are conditioning
on $\mc{G}_i$ (interpreted for the current application of Lemma \ref{de})
and $P_I$ by definition of the event $\mc{H}_i$.
Thus we estimate $\mb{E}(Y^-_{1,I}|\mc{G}_i \wedge \mc{H}_i)$ as
$$((1 \pm e(t)/s_e)(q(t) \pm \tT(t)/s_e)n^2)^{-1} \cdot
 (1 \pm e(t)/s_1)(q(t) \pm \tT(t)/s_1)\aA^2$$
$$ \cdot(1 \pm 2e(t)/s_e)(a_H (2t)^{e_H-2}q(t) \pm \tT(t)/s_e \pm n^{-5\eps}) p^{-1}.$$
Recalling that $s=pn^2$, $y^-_q(t) = a_H (2t)^{e_H-2}q(t)$ and $h_q(t)=(eq+\gamma)'(t)$
we see that we have the required condition
$$\mb{E}(Y^-_{1,I}|\mc{G}_i \wedge \mc{H}_i) = (y^-_q(t) \pm h_q(t)/4s_1)\aA^2/s.$$

The boundedness hypothesis follows immediately from the 
definition of $I$-open pairs.
Note that we can arrange
for \( s_1^2 k_1 n^\eps = \alpha n^{5 \eps} < n \), since \( \eps \) is small.
The remaining conditions of Lemma \ref{de} follow by similar calculations
as in the proof of Theorem \ref{track}. \qed

Next we show that a similar estimate holds for the number of open pairs in $I$.

\begin{lemma} \label{open-I}
If $H$ has smooth independence, then with high probability, for every set $I$ of size $\aA$,
the number of open ordered pairs in $I$ at step $i$
is $Q_I(i) = (1 \pm e(t)n^{-2\eps})(q(t) \pm 2n^{-2\eps})\aA^2$.
\end{lemma}

\nib{Proof.}
We need to estimate the number of ordered pairs in $I$ that are $I$-open but not open.
By Corollary \ref{closed-estimate} we can bound the number of pairs closed by
any edge by $p^{-1}\log n$ (say). By smooth independence we can assume that
$|P_I| < n^{-5\eps}p^{-1}$, so at most $n^{-5\eps}p^{-1} \cdot p^{-1}\log n$ pairs
in $I$ are closed but $I$-open. The required bound follows from these
estimates and Lemma \ref{I-open}. \qed

Finally, we can show that the independence number of the process at time $m$ is at most $\aA$.

\begin{lemma} \label{indep}
If $H$ has smooth independence, then with high probability,
at time $m$ every set $I$ of size $\aA$ contains at least one edge.
\end{lemma}

\nib{Proof.}
At step $i+1$ the process chooses an edge uniformly at random from one of
the $Q(i)$ open ordered pairs. Since $Q_I(i)$ of these belong to $I$,
it fails to choose an edge in $I$ with probability $1 - Q_I(i)/Q(i)$.
Multiplying these probabilities and taking a union bound over $I$ we can
bound the probability that there is an independent set $I$ of size $\aA$
by $p_\aA = \binom{n}{\aA} \max_I \prod_{i=1}^m (1 - Q_I(i)/Q(i))$.
By Theorem \ref{track} and Lemma \ref{open-I} we have
$$Q_I(i)/Q(i) = ((1 \pm e(t)/s_e)(q(t) \pm \tT(t)/s_e)n^2)^{-1}
(1 \pm e(t)n^{-2\eps})(q(t) \pm 2n^{-2\eps})\aA^2.$$
Recalling that $s_e = n^{1/2e_H-\eps}$ and $\mu$ is chosen small enough that
$q(t)^{-1}$ and $e(t)$ are at most $n^\eps$ for $t \le t_{\max}$ we can
estimate $Q_I(i)/Q(i) = (1 \pm 10n^{-\eps})(\aA/n)^2$.
Therefore
$$\log p_\aA = \aA(\log n - \log \aA + 1 + O(1/n)) - m\big((1 \pm 10n^{-\eps})(\aA/n)^2 \pm 2(\aA/n)^4\big).$$
Also, since $m = \mu (\log n)^{1/(e_H-1)} pn^2$ and $\aA = 3\mu^{-1} (\log n)^{1-1/(e_H-1)} p^{-1}$
we have $m(\aA/n)^2 = 3\aA\log n$. 
Thus we obtain $$\log p_\aA < -\aA\log n = -3\mu^{-1}(\log n)^{2-1/(e_H-1)}p^{-1},$$
so $p_\aA < \exp (-n^{1/e_H})$ (say), as required. \qed

\section{Independence number and Ramsey bounds}

In this section we show that cliques and cycles both have the smooth independence property.
By Lemma \ref{indep}, this is enough to prove Theorems \ref{k-indep} and \ref{cycle-indep},
and then Theorem \ref{ramsey} follows immediately from Theorem \ref{k-indep}.
We will also show that a graph $H$ satisfying the hypothesis of Theorem \ref{balanced-indep}
has smooth independence, which is enough to prove that theorem.

We start with cycles, where we deduce smooth independence from a path-counting argument.

\begin{lemma} \label{cycle-smooth}
The $\ell$-cycle $C_\ell$ has smooth independence for $\ell \ge 4$.
\end{lemma}

\nib{Proof.}
Suppose $I \sub [n]$ is a set of $\aA$ vertices and let $P_I \sub E(m)$
be the ordered edges at time $t_{\max}$ that are $I$-bad.
We need to show that with high probability $|P_I| < n^{-5\eps}p^{-1}$ for all such $I$.
Consider the contrary event that $|P_I| \ge n^{-5\eps}p^{-1}$, i.e. there are at least $n^{-5\eps}p^{-1}$
ordered edges that each close at least $n^{-5\eps}p^{-1}$ ordered pairs in $I$.
Then there is some ordered pair of edges $uv$, $xy$ of $C_\ell$
and $P'_I \sub P_I$ with $|P'_I| \ge \ell^{-1}n^{-5\eps}p^{-1}$ such that for every edge $cd$ in $P'_I$
there are at least $\ell^{-1}n^{-5\eps}p^{-1}$ embeddings $f$ of $C_\ell \sm uv$
with $f(x)=c$, $f(y)=d$ and $f(u), f(v) \in I$.

Set $I_0=I$ and for $1 \le j \le \ell-2$ define
$$I_j = \{v: |N_{G(m)}(v) \cap I_{j-1}| > n^{-10\eps} pn\}.$$
By Theorem \ref{track} the degree of any vertex at time $t$ is
$(1 \pm e(t)/s_e)(2t \pm 1/s_e)pn$.
Now $p = n^{-(\ell-2)/(\ell-1)}$ and $t \le t_{\max} = \mu(\log n)^{1/(\ell-1)}$,
so $pn = n^{1/(\ell-1)}$ and we can bound all degrees by $(n\log n)^{1/(\ell-1)}$.
It follows that there are at most $(n\log n)^{j/(\ell-1)}$ paths
of length $j$ starting at any given vertex, for any $j$.
Also, if $v \notin I_j$ we can improve on this estimate when counting
paths of length $j$ that start at $v$ and end in $I$.
To see this, consider choosing the vertex sequence of such a path starting at $v$,
say $v = v_{j-1}, \cdots, v_0 \in I$. At each step we have at most $(n\log n)^{1/(\ell-1)}$
choices, and there must be some $j-1 \ge j' \ge 1$ where
$v_{j'} \notin I_{j'}$ but $v_{j'-1} \in I_{j'-1}$,
when by definition we have at most $n^{-10\eps} pn$ choices.
This gives at most $j(n^{-10\eps} pn)((n\log n)^{(j-1)/(\ell-1)}) < n^{-9\eps} n^{j/(\ell-1)}$
paths of length $j$ that start at $v$ and end in $I$.

Suppose without loss of generality that removing $uv$ and $xy$ from the cycle
leaves a path of length $\ell_1$ joining $u$ to $x$ and a path of length $\ell_2$
joining $v$ to $y$, with $\ell_1+\ell_2=\ell-2$ and $\ell_1>0$
(we might have $\ell_2=0$, i.e. $v=y$).
We claim that for any edge $cd$ in $P'_I$ we must have
$c \in I_{\ell_1}$ and $d \in I_{\ell_2}$. For suppose that $c \notin I_{\ell_1}$.
Then there are at most $n^{-9\eps} n^{\ell_1/(\ell-1)}$
paths of length $\ell_1$ that start at $c$ and end in $I$.
Also, there are at most $(n\log n)^{\ell_2/(\ell-1)}$ paths of length $\ell_2$
that start at $d$ and end in $I$. Thus we bound the number of embeddings
$f$ of $C_\ell \sm uv$ with $f(x)=c$, $f(y)=d$ and $f(u), f(v) \in I$
by $n^{-9\eps} n^{\ell_1/(\ell-1)} \cdot (n\log n)^{\ell_2/(\ell-1)} < \ell^{-1}n^{-5\eps}p^{-1}$,
contradiction. Thus we have $c \in I_{\ell_1}$, and the same argument gives $d \in I_{\ell_2}$.

Now by Lemma \ref{degree}, with high probability we have
$|I_j| \le \aA (8^{-1}\eps n^{-10\eps} pn)^{-j} < n^{1-(j+1)/(\ell-1)+11j\eps}$
for $1 \le j \le \ell-2$ and every $I$ of size $\aA$.
Then by Lemma \ref{edges}, with high probability we have
$$e(I_{\ell_1},I_{\ell_2}) <  \max\{4\eps^{-1}(|I_{\ell_1}|+|I_{\ell_2}|),p|I_{\ell_1}||I_{\ell_2}|n^{2\eps}\}.$$
This is less than $n^{-1/\ell}p^{-1}$ unless $\ell_2=0$.
Also, if $\ell_2=0$ then $\ell_1=\ell-2$, so $|I_{\ell_1}| < n^{11\ell\eps}$
and we can bound the number of edges incident to $I_{\ell_1}$
by $|I_{\ell_1}|(n\log n)^{1/(\ell-1)} < n^{1/(\ell-1) + 12\ell\eps}$.
Either way we have $e(I_{\ell_1},I_{\ell_2}) < \ell^{-1}n^{-5\eps}p^{-1} \le |P'_I|$,
by our earlier assumption, which contradicts the fact any edge $cd$ in $P'_I$ 
has $c \in I_{\ell_1}$ and $d \in I_{\ell_2}$. Therefore with high probability we have
$|P_I| < n^{-5\eps}p^{-1}$ for all $I$, i.e. $H$ has the smooth independence property. \qed

For cliques, we first consider the case $H=K_s$ for some $s \ge 6$.
Then $p = n^{-2/(s+1)}$. Consider any two edges $uv$, $xy$ of $H$ and let $H^- = H \sm uv$.
We have $S_{xy,H^-} = p^{-1}$ and for $s \ge 6$ we have $S_{xy,H^-[B]} \ge p^2 n > p^{-1}$
for any $B$ with $xy \subn B \subn V_H$, i.e. $(xy,H^-)$ is strictly balanced.
We show that this more general property suffices for smooth independence. Note that if
$H$ is any graph such that $(xy, H^-)$ is strictly balanced for all $xy,uv \in E_H$ 
then \(H\) has minimum degree at least 3. (To see this, assume for a contradiction that $ d_H(u) = 2 $ and 
consider an extension $(xy, H^-)$ where $ u \not\in xy $.)

\begin{lemma} \label{balanced-smooth}
Suppose that $(xy, H \sm uv)$ is strictly balanced
for any two edges $uv$, $xy$ of $H$. Then $H$ has smooth independence.
\end{lemma}

\nib{Proof.}
Suppose $I \sub [n]$ is a set of $\aA$ vertices and let $P_I \sub E(m)$
be the ordered edges at time $t_{\max}$ that are $I$-bad.
We need to show that with high probability $|P_I| < n^{-5\eps}p^{-1}$ for all such $I$.
Consider the contrary event that $|P_I| \ge n^{-5\eps}p^{-1}$, i.e. there are at least $n^{-5\eps}p^{-1}$
ordered edges that each close at least $n^{-5\eps}p^{-1}$ ordered pairs in $I$.
Then there is some ordered pair of ordered edges $uv$, $xy$ of $H$ with $u \notin \{x,y\}$
and $P'_I \sub P_I$ with $|P'_I| \ge (2e_H)^{-1}n^{-5\eps}p^{-1}$ such that for every edge $cd$ in $P'_I$
there are at least $(2e_H)^{-1}n^{-5\eps}p^{-1}$ embeddings $f$ of $H^- = H \sm uv$
with $f(x)=c$, $f(y)=d$ and $f(u), f(v) \in I$.

Write $H^- = H \sm uv$. Since $(xy,H^-)$ is strictly balanced we have $S_{B,H^-} < 1$ for any $B$
with $xyu \sub B \subn V_H$. Applying Lemma \ref{extend-general}, we see that
for any $a,c,d \in [n]$ there are at most $n^{4e_H\eps}$ embeddings $f$ of
$H^- = H \sm uv$ with $f(x)=c$, $f(y)=d$ and $f(u) = a$.
For each edge $cd \in P'_I$ let $U_{cd}$ be the set of vertices $a \in I$ such that there is
at least one embedding $f$ of $H^- = H \sm uv$ with $f(x)=c$, $f(y)=d$ and $f(u) = a$.
By definition of $P'_I$ we must have
$$|U_{cd}| > (2e_H)^{-1}n^{-5\eps}p^{-1}/n^{4e_H \eps} > n^{-10e_H\eps}p^{-1}$$
(say) for every edge $cd \in P'_I$. Next we need the following claim.

\nib{Claim.} $|U_{cd} \cap U_{c'd'}| < n^{-1/e_H}p^{-1}$ for any two edges $cd, c'd' \in P'_I$.

\nib{Proof.} Consider two embeddings $f_1, f_2$ of $H^-$ such that $f_1(x)=c$, $f_1(y)=d$,
$f_2(x)=c'$, $f_2(y)=d'$ and $f_1(u)=f_2(u)=a$.
Let $C' = f_1(V_H) \cap f_2(V_H)$ and $H' = H \sm \{uv,xy\}$.
Let $W$ be the join of $J_1=H'$ and $J_2=H'$ formed by identifying the sets
$f_1^{-1}(C')$ and $f_2^{-1}(C')$ as a single set $C$ on which $f_1$ and $f_2$ agree.
Note that we have $u \in C$.
For ease of notation we let $x,y$ denote the copies of $x,y$ in $J_1$
and $x',y'$ the copies of $x,y$ in $J_2$.
Let $A = \{x,y\} \cup \{x',y'\}$. Since $cd \ne c'd'$ we have $|A| \ge 3$.
Define $\phi:A \to [n]$ by $\phi(x)=c$, $\phi(y)=d$, $\phi(x')=c'$, $\phi(y')=d'$.
We want to estimate $N_{\phi,W}$.
The argument is very similar to that in Lemma \ref{closure}.
Choose $B$ with $A \sub B \sub V_W$ maximising $S_{B,W}$.
We have cases depending on how $V_{J_1}$ and $V_{J_2}$ intersect.
If $f_1(V_H) = f_2(V_H)$, i.e. $V_{J_1}=V_{J_2}$,
then we have $S_{B,W} \le S_{B,H^-} \le 1 < p^{-1}$,
since $(xy,H^-)$ is strictly balanced and $|A| \ge 3$.
We henceforth suppose that $f_1(V_H) \ne f_2(V_H)$.
Define $B_1 = B \cap V_{J_1}$ and $B_2 = B \cap V_{J_2}$.
Next we consider the case \( V_{J_1} \subseteq V_{J_2} \cup A \).
If \( B_2 \neq \{x',y'\} \) then we have \( S_{B,W} \le S_{B_2,J_2} \le 1 < 1/p \)
because \( (x'y', J_2) \) is strictly balanced.  If \( B_2 = \{x',y'\} \)
then we note that, since \( H \) has minimum degree at least 3 and $V_{J_1} \sm V_{J_2} \ne \emptyset$, 
we have \( S_{B,W} \le p S_{B_2,J_2} \le 1 < 1/p \).  The analogous argument
handles the case \( V_{J_2} \subseteq V_{J_1} \cup A \).

Now suppose that \(V_{J_1} \setminus ( V_{J_2} \cup A) \)
and \(V_{J_2} \sm ( V_{J_1} \cup A) \) are non-empty.
We consider subcases according to \( B_1 \) and \(B_2\).
The first subcase is $B_1 \cup C \ne V_{J_1}$.
Then we have $S_{B_1 \cup C,J_1} = S_{B_1 \cup C,H^-} < 1$,
since $u \in C$ and $(xy,H^-)$ is strictly balanced.
Also $S_{B_2,J_2} \le S_{x'y',H^-} = p^{-1}$,
so $S_{B,W} = S_{B_2,J_2} S_{B_1 \cup C,J_1} < p^{-1}$.
The second subcase is $B_2 \cup C \ne V_{J_2}$,
when a similar argument gives
$S_{B,W} = S_{B_1,J_1} S_{B_2 \cup C,J_2} < p^{-1}$.
Finally, the third subcase is $B_1 \cup C = V_{J_1}$ and $B_2 \cup C = V_{J_2}$.
Then $B_1$ contains $V_{J_1} \sm (A_1 \cup C)$ and
$B_2$ contains $V_{J_2} \sm (A_2 \cup C)$,
which are both non-empty.
Since $(xy,H^-)$ is strictly balanced
we have $S_{B_1,J_1} \le 1$ and $S_{B_2 \cup C,J_2} \le 1$,
and so $S_{B,W} = S_{B_1,J_1} S_{B_2 \cup C,J_2} \le 1$.
In all cases we have $S_{B,W} < p^{-1}$, so $S_{B,W} \le n^{-1/(e_H-1)}p^{-1}$,
since it is an integer power of $n^{-1/(e_H-1)}$. Now Lemma \ref{extend-general}
gives $N_{\phi,W} < n^{4e_W \eps -1/(e_H-1)}p^{-1}$.
Summing over all possible joins $W$ we estimate
$|U_{cd} \cap U_{c'd'}| < n^{-1/e_H}p^{-1}$, which proves the claim. \qed

Returning to the proof of the lemma,
we now set $\omega = n^{11e_H\eps}$ and choose $\omega$ edges of $P'_I$,
say $c_1d_1,\cdots,c_{\omega}d_{\omega}$. Recall that
$|U_{cd}| > n^{-10e_H\eps}p^{-1}$ for every $cd \in P'_I$.
Then $|U_{c_id_i} \sm \cup_{j<i} U_{c_jd_j}| > n^{-10e_H\eps}p^{-1} - in^{-1/e_H}p^{-1}$
for $1 \le i \le \omega$ by the claim. This gives
$$|\cup_{i=1}^\omega U_{c_id_i}| > \omega n^{-10e_H\eps}p^{-1} - \frac{1}{2}\omega^2 n^{-1/e_H}p^{-1}
> n^{\eps} p^{-1},$$ 
say. But by definition the sets
$U_{c_id_i}$ are contained in $I$, 
for which $|I|=\aA = 3 \mu^{-1} (\log n)^{1-1/(e_H-1)} p^{-1}$ is too small.
This contradiction shows that we cannot have $|P_I| \ge n^{-5\eps}p^{-1}$ for some $I$
holding together with the bounds used from Lemma \ref{extend-general}.
These bounds hold with high probability, so with high probability we have
$|P_I| < n^{-5\eps}p^{-1}$ for all $I$, i.e. $H$ has the smooth independence property. \qed

The two arguments above can be generalised to prove smooth independence for a wider
class of graphs $H$. However, for the sake of brevity and clarity, 
we restrict our attention to these simple cases here.
We complete the discussion of cliques by showing that $K_5$ has smooth independence.
(The independence numbers for the $K_3$-free and $K_4$-free processes have
already been obtained in \cite{B}.)

\begin{lemma} \label{k5-smooth}
$K_5$ has smooth independence.
\end{lemma}

\nib{Proof.}
Write $H=K_5$. We argue as in the proof of Lemma \ref{balanced-smooth}.
Consider $I$, $P_I$, $uv$, $xy$, $P'_I$, $U_{cd}$ as defined in that proof.
Now $(xy,H^-)$ is not strictly balanced, but do we have $S_{B,H^-} \le 1$
for any $B$ with $xyu \sub B \sub V_H$, so for every $cd \in P'_I$
we still obtain the bound $|U_{cd}| > n^{-10e_H\eps}p^{-1}$. Following that proof,
our next step is to show that $|U_{cd} \cap U_{c'd'}| < n^{-1/e_H}p^{-1}$
for any two edges $cd, c'd' \in P'_I$. In fact we will obtain a much stronger bound.
Consider two embeddings $f_1, f_2$ of $H$ such that $f_1(x)=c$, $f_1(y)=d$,
$f_2(x)=c'$, $f_2(y)=d'$ and $f_1(u)=f_2(u)=a$.
Define $C'$, $H'$, $J_1$, $J_2$, $W$, $x'$, $y'$, $A$ and $\phi$ as before.
Choose $B$ with $A \sub B \sub V_W$ maximising $S_{B,W}$.
Note that for any $K$ with $xy \subn K \sub V_H$ we have $S_{K,H^-} \le 1$.
So if $V_{J_1}=V_{J_2}$ we have $S_{B,W} \le 1$.
Otherwise we consider cases according to $B_1 = B \cap V_{J_1}$ and $B_2 = B \cap V_{J_2}$.
Since $S_{B,W} = S_{B_1,J_1} S_{B_2 \cup C,J_2}$, $S_{B,W} = S_{B_2,J_2} S_{B_1 \cup C,J_1}$
and $|B_1 \cup C|, |B_2 \cup C| \ge 3$ we see that $S_{B,W} \le 1$,
except possibly in the case $B_1 = \{x,y\}$ and $B_2=\{x',y'\}$.
In this case we note that there is an edge from $a \in C$ to $B_2$ that
is not contained in $J_1$, so $S_{B,W} \le pS_{B_1,J_1} = 1$.
In all cases we have $S_{B,W} \le 1$ and so $N_{\phi,W} < n^{4e_W}$.
Summing over all possible joins $W$ we estimate $|U_{cd} \cap U_{c'd'}| < n^{5e_W}$, say.
Now the remainder of the proof follows as in Lemma \ref{balanced-smooth}. \qed

\section{Concluding remarks}

We have restricted our attention in this paper to those aspects
of the $H$-free process needed for our applications
to Ramsey and Tur\'an bounds. However, we also view this work as
the first stage in the study of this process as a model of independent
interest. In the course of our arguments we have already described some
properties of the model via our asymptotic formulae for trackable
extension variables; for example, we have shown that for fixed
graphs \( \GG \) that do not contain \(H\) as a subgraph,
excluding `critical' cases, the number of copies of 
\( \GG \) in \(G(i) \) is roughly the same as the number of copies of \(\GG\) in
the unconstrained random graph \( G(n,i) \).  
In principle, one may ask for analogues in the graph \( G(i) \) produced
by the \(H\)-free process of any property known to hold in \( G(n,i) \).
But the most natural next steps are continued investigation of the
independence number and development of upper bounds on the number of
steps in the \(H\)-free process.  For independent sets, there
are other classes of graphs covered by our methods, but for clarity
we have restricted our attention to certain concrete settings
rather than stating a complicated general theorem. One might hope that
any strictly $2$-balanced graph can be analysed by these methods.
With respect to upper bounds, we believe that the number of steps 
in the \(H\)-free process is at most a constant times the lower bound 
we establish here for any strictly 2-balanced \(H\).
In fact, we are even prepared to make this conjecture for the
degree of each vertex.
\begin{conj}
For any strictly 2-balanced graph $H$ there is a constant $C$
so that with high probability the maximal \(H\)-free graph $G$ 
on \(n\) vertices produced by the \(H\)-free process has maximum degree
\[ \Delta(G) < C n^{1 - (v_H-2)/(e_H-1)} (\log n)^{1/(e_H-1)}. \]
\end{conj}
\noindent
For the triangle-free process this follows from the
bound on the independence number (see \cite{B}), but in general it is
a separate question.  The later evolution of the process, 
where Theorem~\ref{track} no longer applies, is also an intriguing
topic for further study.


\begin{thebibliography}{99}

\bibitem{AKS}
M. Ajtai, J. Koml\'os and E. Szemer\'edi, A note on Ramsey numbers,
{\em J. Combin. Theory Ser. A} {\bf 29} (1980), 354--360.

\bibitem{ABK}
N. Alon, S. Ben-Shimon and M. Krivelevich, 
A note on regular Ramsey graphs, arXiv:0812.2386v1. 

\bibitem{ARS}
N. Alon, L. R\'{o}nyai and T. Szab\'{o}, 
Norm-graphs: variations and applications,
{\em J. Combin. Theory Ser. B} {\bf 76} (1999), 280--290. 

\bibitem{AS} N. Alon and J. Spencer,
{\em The probabilistic method}, second edition,
Wiley, New York, 2000.

\bibitem{Ap} T. M. Apostol,
An elementary view of Euler's summation formula,
{\em Amer. Math. Monthly} {\bf 106} (1999), 409--418.

\bibitem{BKR} A. D. Barbour, M. Karo\'nski and A. Ruci\'nski,
A central limit theorem for decomposable random variables, 
with applications to random graphs,
{\em J. Combin. Theory Ser. B} {\bf 47} (1989), 125--145.

\bibitem{B} T. Bohman, The triangle-free process, {\em Advances in Mathematics}
{\bf 221} (2009), 1653-1677.

\bibitem{BR}
B. Bollob\'as and O. Riordan, Constrained graph processes,
{\em Electronic J. Combin.} {\bf 7} (2000) R18.

\bibitem{Br} W. G. Brown,
On graphs that do not contain a Thomsen graph,
{\em Canad. Math. Bull.} {\bf 9} (1966), 281--289.

\bibitem{CLRZ}
Y. Caro, Y. Li, C. C. Rousseau and Y. Zhang,
Asymptotic bounds for some bipartite graph - complete graph Ramsey numbers,
{\em Disc. Math.} {\bf 220} (2000), 51--56.

\bibitem{D}
R. Durrett, {\em Random Graph Dynamics}, Cambridge Univ. Press, 2007.

\bibitem{E} P. Erd\H{o}s,
Graph theory and probability, II,
{\em Canad. J. Math.} {\bf 13} (1961), 346--352.

\bibitem{E2} P. Erd\H{o}s,  
Extremal problems in number theory, combinatorics and geometry,  
{\em Proc. ICM}, PWN, Warsaw, 1984, 51--70.

\bibitem{ESi} P. Erd\H{o}s and M. Simonovits,
Some extremal problems in graph theory,
{\em Coll. Math. Soc. J\'anos Bolyai} {\bf 4} (1969), 377--390.

\bibitem{ESp} P. Erd\H{o}s and J. H. Spencer,
{\em Probabilistic Methods in Combinatorics},
Academic Press, 1974.

\bibitem{ESt} P. Erd\H{o}s and  A.H. Stone,
On the structure of linear graphs,
{\em Bull. Amer. Math. Soc.} {\bf 52} (1946), 1087--1091.

\bibitem{ESW} P. Erd\H{o}s, S. Suen and P. Winkler,
On the size of a random maximal graph,
{\em Random Structures Algorithms} {\bf 6} (1995), 309--318.


\bibitem{F} Z. F\"uredi,
Tur\'an type problems, in: {\em  Surveys in combinatorics},
London Math. Soc. Lecture Note Ser. 166, Cambridge Univ. Press,
Cambridge, 1991, 253--300

\bibitem{F2} Z. F\"uredi,
New asymptotics for bipartite Tur\'an numbers,
{\em J. Combin. Theory Ser. A} {\bf 75} (1996), 141--144. 

\bibitem{F3} Z. F\"uredi,
An upper bound on Zarankiewicz' problem,
{\em Combin. Probab. Comput.} {\bf 5} (1996), 29--33. 

\bibitem{G}
D. Grable, On random greedy triangle packing,
{\em Electronic J. Combin.} {\bf 4} (1997), R11. 

\bibitem{GRS} R. L. Graham, B. L. Rothschild and J. H. Spencer,
{\em Ramsey Theory}, Wiley, New York, 1990.

\bibitem{K} J. H. Kim,
The Ramsey number $R(3,t)$ has order of magnitude $t^2/\log t$,
{\em Random Structures Algorithms} {\bf 7} (1995), 173--207.

\bibitem{KST} T. K\"{o}vari, V. T. S\'{o}s and P. Tur\'{a}n, 
On a problem of K. Zarankiewicz, {\em Colloquium Math.} 
{\bf 3} (1954), 50--57.

\bibitem{LZ} Y. Li and W. Zang,
The independence number of graphs with a forbidden cycle and Ramsey numbers,
{\em J. Combin. Opt.} {\bf 7} (2003), 353--359. 

\bibitem{OT} D. Osthus and A. Taraz,
Random maximal $H$-free graphs,
{\em Random Structures Algorithms} {\bf 18} (2001), 61--82.

\bibitem{R} A. Ruci\'nski,
Recent developments in random graphs, 
Proceedings of the International Summer School on Probability and Statistics, 
Varna (1994), online at:
{\tt http://www.staff.amu.edu.pl/$\sim$rucinski/papers/43.pdf}
 
\bibitem{RW} A. Ruci\'nski and N. Wormald,
Random graph processes with degree restrictions,
{\em Combin. Probab. Comput.} {\bf 1} (1992), 169--180.

\bibitem{Se1} T. G. Seierstad,
A central limit theorem via differential equations, 
{\em Ann. Appl. Probab.} {\bf 19} (2009), 661–675.

\bibitem{Se2} T. G. Seierstad,
Stronger large deviation bounds for Wormald's differential equation method, submitted.

\bibitem{Sp1} J. Spencer,
Asymptotic lower bounds for Ramsey functions,
{\em Disc. Math.} {\bf 20} (1997), 69--76.

\bibitem{Sp2} J. Spencer, Counting extensions,
{\em J. Combin. Theory Ser. A} {\bf 55} (1990), 247--255.

\bibitem{Su} B. Sudakov,
A note on odd cycle-complete graph Ramsey numbers, 
{\em Electronic J. Combin.} {\bf 9} (2002), N1.

\bibitem{T} P. Tur\'an,  Eine Extremalaufgabe aus der
Graphentheorie, {\em Mat. Fiz. Lapok} {\bf 48} (1941), 436--452.

\bibitem{Wz} G. Wolfovitz,
Lower bounds for the size of random maximal $H$-free graphs,
{\em Electronic J. Combin.} {\bf 16} (2009), R4.
  
\bibitem{W} N. C. Wormald,
The differential equation method for random graph processes
and greedy algorithms, in:
{\em Lectures on Approximation and Randomized Algorithms},
PWN, Warsaw, 1999, 73--155.

\end{thebibliography}
\end{document}